\documentclass{amsart}
\usepackage{amsmath}
  \usepackage{paralist}
  \usepackage{graphics} 
  \usepackage{epsfig} 
\usepackage{graphicx}  \usepackage{epstopdf}
 \usepackage[colorlinks=true]{hyperref}
 \usepackage{bbm, dsfont, mathrsfs}
 \usepackage{verbatim}
\hypersetup{urlcolor=blue, citecolor=red}
\DeclareMathOperator*{\esssup}{ess\,sup}

\newtheorem{theorem}{Theorem}[section]
\newtheorem{corollary}{Corollary}

\newtheorem{lemma}[theorem]{Lemma}
\newtheorem{proposition}{Proposition}

\theoremstyle{definition}
\newtheorem{definition}[theorem]{Definition}
\newtheorem{remark}{Remark}

\newcommand{\eps}[1]{{#1}_{\varepsilon}}

\title[Interface Perturbations and non-linearity]
      {Notes on the nonlinear dependence of a multiscale coupled system \\
      with respect to the interface}

\author[Fernando A Morales]{}


 \email{famoralesj@unal.edu.co}
 \email{show@math.oregonstate.edu}

\thanks{The author was supported by the project HERMES 17194 from Universidad Nacional de Colombia}

\begin{document}

\def\A{{\mathcal A}}
\def\B{{\mathcal B}}
\def\C{{\mathcal C}}
\def\E{{\mathcal E}}
\def\F{{\mathcal F}}
\def\Q{{\mathcal Q}}
\def\V{\mathbf V}
\def\W{\mathbf W}
\def\H{\mathbf H}
\def\X{\mathbf X}
\def\Y{\mathbf Y}
\def\0{\mathbf 0}
\def\f{\mathbf f}
\def\g{\mathbf g}
\def\n{\boldsymbol {\widehat{n } } }
\def\q{\mathbf q}
\def\r{\mathbf r}
\def\s{\mathbf s}
\def\u{\mathbf u}
\def\v{\mathbf v}
\def\w{\mathbf w}
\def\x{\mathbf x}
\def\y{\mathbf y}
\def\div{\boldsymbol{\nabla} \cdot}
\def\grad{\boldsymbol{\nabla}}
\def\ygrad{\boldsymbol{\nabla_y}}
\def\del{\boldsymbol{\partial}}
\def\eps{\varepsilon}
\def\cC{\mathbb C}
\def\hH{\mathbb H}
\def\vV{\mathbb V}
\def\:{\negthinspace : \negthinspace}
\def\bone{{ 1 }}
\def\chix{{\raise.5ex\hbox{$\chi$}}}
\def\myskip{\noalign{\vskip6pt}}
\def\norm#1{\| #1 \|}
\def\vector#1{\boldsymbol{ #1 }}
\def\laplace#1{\overline{ #1 }}
\def\Re{\mbox{\rm I\negthinspace R}}
\def\Co{\mbox{C\kern-.47em\vrule height1.5ex width.2ex}}
\def\wconv{\overset{w\,\,}\rightharpoonup}
\def\sconv{\rightarrow}
\def\Dom{\operatorname{Dom}}
\def\Rg{\operatorname{Rg}}
\def\Ker{\operatorname{Ker}}
\def\gradth{\boldsymbol{\widetilde{\nabla} } }

\def\ueps{u}
\def\uepsz{u^{\zeta}}
\def\uepsp{u^{\,\epsilon '}}
\def\uepsone{u^{1}}
\def\uepszone{u^{1}_{\zeta}}
\def\uepstwo{u^{2}}
\def\uepsztwo{u^{2}_{\zeta}}
\def\etaeps{\eta^{\,\epsilon}}
\def\mueps{\mu^{\,\epsilon}}
\def\muepsone{\mu^{\,\epsilon,\,1}}
\def\muepstwo{\mu^{\,\epsilon,\,2}}
\def\etaueps{\eta^{\,\epsilon}}
\def\N{\boldsymbol {\mathbbm{N} } }
\def\R{\boldsymbol{\mathbbm {R} } }
\def\eversor{\mathbf{\widehat{e}}}

\def\xthilde{\widetilde{\mathbf{x}}}
\def\Xthilde{\widetilde{\mathbf{X}}}
\def\xn{x_{\scriptscriptstyle N}}
\def\pert{\mathscr{ T } }
\def\apert{\mathscr{ T } }
\def\defining{\overset{\mathbf{def}}=}
\def\Omrz{\Omega_{\,1}^{\,\zeta}}
\def\Omfz{\Omega_{\,2}^{\,\zeta}}
\def\Gammaz{\Gamma^{\,\zeta}}
\def\Gammazp{\Gamma^{\,\zeta}_{+}}
\def\Gammazm{\Gamma^{\,\zeta}_{-}}
\def\Gammap{\Gamma_{+}}
\def\Gammam{\Gamma_{-}}
\def\Uz{\mathcal{U}^{\,\zeta}}
\def\Oz{\mathcal{O}^{\zeta}}

\def\sOmrz{\Omega_{1}^{\zeta}}
\def\sOmfz{\Omega_{2}^{\zeta}}
\def\sGammaz{\Gamma^{\zeta}}
\def\sGammazp{\Gamma^{\zeta}_{+}}
\def\sGammazm{\Gamma^{\zeta}_{-}}
\def\sUz{\mathcal{U}^{\zeta}}

\def\Oone{\mathcal{O}_{1}}
\def\Otwo{\mathcal{O}_{2}}
\def\Othree{\mathcal{O}_{3}}

\def\deps{\mathcal{D}^{\epsilon}}
\def\w{\mathbf{w}}
\def\wtan{\mathbf{w}_{T}}
\def\wnorm{\mathbf{w}_{N}}

\def\Hdiv{\mathbf{H}_{\,\text{div}}}
\def\Ltwo{\mathbf{L}^{2}}
\def\Vz{\mathbf{V}^{\zeta}}
\def\Qz{Q^{\,\zeta}}

\def\vone{\textbf{v}^{1}}
\def\uone{\textbf{u}^{1}}
\def\wone{\textbf{w}^{\,1}}
\def\vtwo{\textbf{v}^{\,2}}
\def\utwo{\textbf{u}^{\,2}}
\def\wtwo{\textbf{w}^{\,2}}
\def\vthree{\textbf{v}^{3}}
\def\uthree{\textbf{u}^{3}}
\def\wthree{\textbf{w}^{\,3}}


\def\pone{p_{ 1}}
\def\qone{q_{ 1}}
\def\rone{r_{ 1}}
\def\ptwo{p_{ 2}}
\def\qtwo{q_{ 2}}
\def\rtwo{r_{ 2}}
\def\pthree{p_{ 3}}
\def\qthree{q_{ 3}}
\def\rthree{r_{ 3}}


\def\Lap{\boldsymbol{\Delta}}
\def\enversor{\boldsymbol {\widehat{e}}_{\scriptscriptstyle N}}
\def\xiv{\boldsymbol {\xi}}
\def\etav{\boldsymbol {\eta}}
\def\ind{\boldsymbol {\mathbbm {1} } }

\maketitle

\centerline{\scshape Fernando A Morales }
\medskip
{\footnotesize
 \centerline{Escuela de Matem\'aticas, Universidad Nacional de Colombia, Sede Medell\'in.}
   \centerline{ Calle 59 A No 63-20, Of 43-106. Medell\'in, Colombia.}
} 

\bigskip


\begin{abstract}
This work studies the dependence of the solution with respect to interface geometric perturbations in a multiscaled coupled Darcy flow system in direct variational formulation. A set of admissible perturbation functions and a sense of convergence are presented, as well as sufficient conditions on the forcing terms, in order to conclude strong convergence statements. For the rate of convergence of the solutions we start solving completely the one dimensional case using orthogonal decompositions on appropriate subspaces. Finally, the rate of convergence question is analyzed in a simple multiple dimensional setting, studying the nonlinear operators introduced by the geometric perturbations. 
\end{abstract}

\section{Introduction}   \label{Sec Introduction}
\begin{figure}[t]
\caption[1]{Original Domain and a Perturbation}\label{Fig Artificial Perturbation}
\centerline{\resizebox{8cm}{7cm}
{\includegraphics{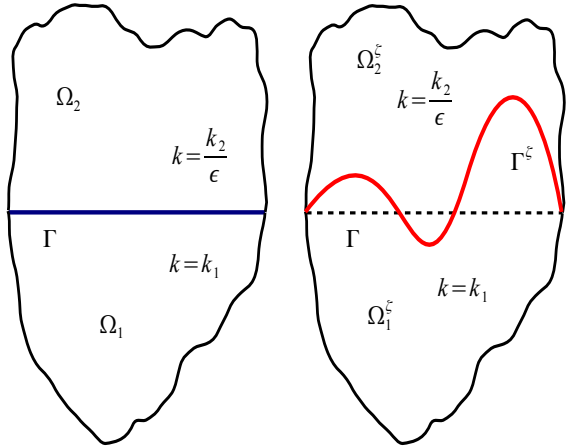} } }
\end{figure}

The study of saturated flow in geological porous media frequently presents natural structures with a dense network of fissures nested in the rock matrix \cite{Show97}. It is also frequent to observe vuggy porous media \cite{ArbLehr2006}, which have the presence of cavities in the rock matrix significantly larger than the average pore size of the medium. This is a multiple scale physics phenomenon, because there are regions of the medium where the flow velocity is significantly larger than the velocity on the other ones. The modeling of the interface between regions is subject to very active research: first, fluid transmission conditions  across the interface are of great importance, see \cite{Saffman, BeaversJoseph} for discussion of the governing laws; see \cite{Gatica2009, ArbBrunson2007} for a numerical point of view; see \cite{ArbLehr2006, Yotov, Gunzburger2010, MoralesShow2} for the analytic approach and  \cite{SP80, CannonMeyer71} for a more general perspective. Second, the placement of the interface is debatable since a boundary layer phenomenon between regions occurs, see \cite{HigashinoStefan, PasterDagan} for discussion. The interface couples regions of slow velocity (order $\mathcal{O}(1)$) and fast flow (order  $\mathcal{O}(1/ \epsilon)$), see figure \ref{Fig Artificial Perturbation}. Hence, its placement and geometric description become an important issue because perturbations of the interface are inevitable. On one hand the geological strata data available are always limited, on the other hand the numerical implementation of models involving curvy interfaces, in most of the cases can only approximate the real surface. Finally, on a very different line, in the analysis of saturated flow in deformable porous media, one of the aspects is the understanding the geometric perturbations of an \emph{interface of reference}. 

Clearly, the continuity of the solution with respect to the geometry of the interface is an important issue which has received very little treatment and mainly limited to flat interfaces. Most of the theoretical achievements in the field of multiscale coupled systems, concentrate their efforts in removing the singularities introduced by the scales using homogenization processes. These techniques can be either formal \cite{Morales1, Levy83}, analytic \cite{Hornung} or numerical \cite{JaffRob05}.

Here, we model the stationary problem with a coupled system of partial differential equations of Darcy flow in both regions, in direct variational formulation.  We simulate the region of fast flow scaling by $\frac{1}{\epsilon}$ the ratio of permeability over viscosity, as in figure \ref{Fig Artificial Perturbation}. It will be assumed that the real interface $\Gamma$ is horizontal flat and the perturbed one $\Gamma^{\zeta}$ is curved. Of course a flat surface will not be perturbed when discretized and seems unrealistic to consider perturbations of it. Our choice is motivated by two reasons: first, for the sake of clarity in the notation, calculation and interpretation of the results. Second, when studying the phenomenon of saturated flow in \emph{deformable porous media} perturbations of a flat surface are of interest. Finally it is important to highlight that the mathematical essentials of the problem are captured in this framework. The paper starts proving in section 2, that the solutions depend continuously with respect to the interface, then it moves to the much deeper question of exploring the dependence itself. In section 3 the one dimensional case the rate of convergence question is solved completely using orthogonal decomposition of adequate subspaces, finally section 4 reveals the highly nonlinear dependence of the solutions with respect to the interface.  
%
%
%
%

We close this section introducing the notation. Vectors are denoted by boldface letters, as are vector-valued functions and corresponding function spaces. We use $\xthilde$ to indicate a vector in $\R^{\! N-1}$; if $\x\in \R^{\! N}$ then the $\R^{\! N-1}\times\{0\}$ projection is identified with $\xthilde\defining(x_{1}, x_{2}, \ldots, x_{\scriptscriptstyle N-1})$ so that $\x = (\xthilde, \xn)$. The symbol $\gradth$ represents the gradient in the first $N-1$ derivatives. Given a function $f: \R^{\! N}\rightarrow \R$ then  $\int_{M} f\,dS$ is the notation for its surface integral on the $\R^{\! N-1}$ manifold $M\subseteq \R^{\! N}$. $\int_{A} f\, d\x$ stands for the volume integral in the set $A\subseteq \R^{\! N}$; whenever the context is clear we simply write $\int_{A} f$. The notation $\Vert\cdot \Vert_{0, A}$, $\Vert \cdot \Vert_{1, A}$ respectively denote the $L^{2}(A)$ and $H^{1}(A)$ norms on the domain $A\subseteq \R^{N}$. $\ind_{A}$ stands for the indicator function of any given set $A$. The $\eversor_{\ell}$ indicates the unitary vector in the $\ell$-th direction for $1\leq \ell\leq N$. We denote by $\boldsymbol{\widehat{\nu} }$ the outwards normal vector to smooth domain in $\R^{\!N}$ and $\n$ indicates the upwards normal vector i.e. $\n \cdot \enversor \geq 0$. Finally, the Lebesgue measure in $\R^{\! N}$ is denoted by $\lambda_{\scriptscriptstyle N}$. 
%
%
%
%
%
%
%
%
%
\section{Formulation and Convergence}\label{Sec numerical stability}
%
%
%
%
%
%
%
%
%
\subsection{Geometric Setting}\label{Sec Geometric Setting}
In the following $\Gamma$ denotes a connected set in $\R^{\! N-1}\times\{0\}$ whose projection onto $\R^{\! N-1}$ is open; from now on we make no distinction between these two domains. Similarly, $\Omega_{\,1}, \Omega_{\,2}$ denote be smooth bounded open regions in $\R^{\! N}$ separated by $\Gamma$, i.e. $\partial\Omega_{\,1}\cap \partial\Omega_{\,2} = \Gamma$; and such that $\text{sgn} (\x\cdot\enversor) = (-1)^{i}$ for each $\x\in \Omega_{i}$, $i = 1, 2$ (see figure \ref{Fig Artificial Perturbation}). Next, we introduce the admissible perturbations of the interface $\Gamma$. 
\begin{definition}\label{Def perturbations concept}
We say the set $\pert(\Gamma, \Omega)$ of piecewise $C^{1}$ perturbations of the interface $\Gamma$ contained in $\Omega$ is given by
\begin{multline}\label{Def perturbations of interface}
\pert(\Gamma, \Omega)\defining\{\zeta\in C\left(\overline{\Gamma}\,\right):
(\xthilde, \zeta(\xthilde\,))\in \Omega \;\,\forall\,\xthilde\in \Gamma\\
\,,\;\zeta\,\vert_{\partial\Gamma} = 0\,\;
\text{and}\;\zeta\;\text{is a piecesise}\;C^{1}\;\text{function}\,\}.
\end{multline}
The interface associated to $\zeta\in \pert(\Gamma, \Omega)$ is given by the set
\begin{subequations}\label{Def domain associated to perturbation}
\begin{equation}\label{Def interface zeta perturbed}
\Gammaz \defining \{(\xthilde, \zeta(\xthilde\,))\in\R^{\!N}:\xthilde\in \Gamma \}.
\end{equation}
The domains associated to $\zeta\in \pert(\Gamma, \Omega)$ are defined by the sets
\begin{equation}\label{Def matrix zeta perturbed}
\Omrz \defining \left\{(\xthilde, \xn)\in\Omega:\xthilde\in \Gamma,\, \xn< \zeta(\xthilde\,) \right\}
\end{equation}
\begin{equation}\label{Def fissure zeta perturbed}
\Omfz\defining\left\{(\xthilde, \xn)\in\Omega:\xthilde\in \Gamma,\,\zeta(\xthilde\,)<\xn \right\}
\end{equation}
\end{subequations}
\end{definition}
\begin{remark} Observe the following facts
\begin{subequations}
\begin{equation}\label{Eq Gammaz cut}
\partial\Omrz\cap\partial\Omfz = \Gammaz ,
\end{equation}
\begin{equation}\label{Eq reconstruction of Omega}
\Omrz\cup\Gammaz\cup\Omfz = \Omega ,
\end{equation}
\begin{equation}\label{Eq Agreement of boundary 1}
\partial\Omrz - \Gammaz = \partial\Omega_{\,1} - \Gamma , 
\end{equation}
\begin{equation}\label{Eq Agreement of boundary 2}
\partial\Omfz - \Gammaz = \partial\Omega_{\,2} - \Gamma .
\end{equation}
\end{subequations}
%
%
\end{remark}
\begin{definition}\label{Def function space setting}
Define the space
\begin{equation}\label{Def space of pressures direct}
V \defining \{u\in H^{\,1}(\Omega): u\vert_{\,\partial\Omega_{1} - \Gamma} = 0\}
\end{equation}
Endowed with the inner product $\langle \cdot, \cdot\rangle: V\times V\rightarrow \R$
\begin{equation}\label{Def norm space of pressures direct}
\langle u, v\rangle_{\scriptscriptstyle V}\defining \int_{\Omega} \grad u\cdot \grad v ,
\end{equation}
and the norm $\Vert u\Vert_{\scriptscriptstyle V} \defining \sqrt {\langle u, u\rangle_{\scriptscriptstyle V} }$.
\end{definition}
\begin{remark}\label{Rem properties of the function space}
Recall that due to the boundary condition defining the space $V$ and the Poincar\'e inequality the $\Vert\cdot\Vert_{\scriptscriptstyle  V}$-norm is equivalent to the standard $H^{1}$-norm.   
\end{remark}
%
%
%
%
%
%
\subsection{The Problems}\label{Sec the problems}
Consider the strong problem
\begin{subequations}\label{Pblm strong original}
   \begin{equation}\label{Eq strong original}   
      - \div \frac{k_{i}}{\epsilon^{\,i -1}} \grad p_{i} = F \quad \text{in}\; \Omega_{i}, \,\;  i = 1, 2, 
    \end{equation}
    with the interface conditions
    \begin{equation}    
    \pone = \ptwo , \quad k_{1}\,\grad \pone \cdot\n - \frac{k_{2}}{\epsilon} \grad \ptwo\cdot \n = f \quad \text{on} \; \Gamma ,           
    \end{equation}
    and the boundary conditions
    \begin{equation}
      \pone = 0 \quad \text{on}\; \partial \Omega_{1} - \Gamma ,\quad 
      \grad\ptwo\cdot\n = 0 \quad \text{on}\; \partial \Omega_{2} - \Gamma .
    \end{equation}      
\end{subequations}
Now, its perturbation in strong form is given by
\begin{subequations}\label{Pblm strong perturbed}
   \begin{equation}\label{Eq strong perturbed}   
      - \div \frac{k_{i} }{\epsilon^{\,i - 1}}\grad q_{i} = F \quad \text{in}\; \Omega_{i}^{\zeta}, \,\;  i = 1, 2, 
    \end{equation}
    with the interface conditions
    \begin{equation}    
    \qone = \qtwo , \quad k_{1}\,\grad \qone \cdot\n - \frac{k_{2}}{\epsilon} \grad \qtwo\cdot \n = f   
    \quad \text{on} \; \Gammaz ,
    \end{equation}
    and the boundary conditions
    \begin{equation}
      \qone = 0 \quad \text{on}\; \partial \Omega_{1}^{\zeta} - \Gammaz ,\quad 
      \grad\qtwo\cdot\n = 0 \quad \text{on}\; \partial \Omega_{2}^{\zeta} - \Gammaz .
    \end{equation}      
\end{subequations}
Both systems above model stationary Darcy flow, coupling the regions depicted in the left and right hand side of figure \ref{Fig Artificial Perturbation} respectively. The coefficients $k_{1}$, $k_{2}$ indicate the permeability in the corresponding domain; for simplicity they will be omitted in the following. The scaling factor $\frac{1}{\epsilon}$ ensures a much higher velocity $\mathcal{O}(\frac{1}{\epsilon})$ in the upper region with respect to the lower region fluid velocity $\mathcal{O}(1)$. The term $F$ stands for fluid sources and $f$ for a normal flux forcing term on the interface; it is assumed that $f$ is well defined $L^{2}$-function on both manifolds $\Gamma$ and $\Gamma^{\zeta}$. Hence, the weak problems in direct formulation are given by
\begin{subequations}
\begin{equation}\label{Pblm Direct Formulation Original}
p\in V:\quad \int_{\Omega_{1}} \grad p\cdot\grad r
+ \frac{1}{\epsilon} \int_{\Omega_{2}} \grad p\cdot\grad r \\
= \int_{\Omega} F\, r + \int_{\Gamma} f\, r\,d\,S\,,\quad \forall\,r\in V
\end{equation}
\begin{equation}\label{Pblm Direct Formulation Perturbed}
q^{\,\zeta}\in V:\quad \int_{\sOmrz} \grad q^{\,\zeta} \cdot\grad r
+ \frac{1}{\epsilon} \int_{\sOmfz} \grad q^{\,\zeta}\cdot\grad r \\
= \int_{\Omega} F\, r + \int_{\sGammaz} f\, r\,d\,S\,,\quad \forall\,r\in V
\end{equation}
\end{subequations}
\begin{theorem}\label{Th well posedness}
   The problems \eqref{Pblm Direct Formulation Original}, \eqref{Pblm Direct Formulation Perturbed} are well-posed.
   \begin{proof}
      Is is a direct application of Lax-Milgram's lemma and the Poincar\'e inequality, see \cite{MoralesShow1} for details.
   \end{proof}
\end{theorem}
%
%
%
%
%
%
\subsection{A-priori Estimates and Weak Convergence}\label{Sec stability estimates}
In this section under reasonable conditions on the forcing terms and the appropriate type of convergence for the perturbations $\zeta$, a-priori estimates on the solutions of problems \eqref{Pblm Direct Formulation Perturbed} as well as weak convergence statements to the solution of problem \eqref{Pblm Direct Formulation Original} are attained. Test equation \eqref{Pblm Direct Formulation Perturbed} with the solution $q^{\zeta}$, due to the boundary conditions of $V$ and the Poincar\'e constant $C_{\scriptscriptstyle \Omega}$ we get
\begin{multline}\label{Ineq first a-priori estimate}
\frac{1}{1 + C^{2}_{\scriptscriptstyle \Omega}}\,\Vert \,q^{\zeta}\Vert_{1,\,\Omega}^{\,2}
\leq  \Vert \grad q^{\zeta}\Vert_{0, \Omega}^{\,2} 
\leq \Vert \grad q^{\zeta}\Vert_{0,\sOmrz}^{\,2}
+ \frac{1}{\epsilon} \, \Vert \grad q^{\zeta}\Vert_{0, \sOmfz}^{\,2}\\[5pt]
\leq \Vert  F\Vert_{0, \Omega} \, \Vert  q^{\zeta}\Vert_{0, \Omega}
+  \int_{\sGammaz}f\, q^{\zeta}\,dS.
\end{multline}
If a sequence of perturbations $\{\zeta_{n}\}\subseteq \pert(\Gamma, \Omega)$ is to be analyzed, conditions on the type of convergence must be specified. For the perturbations, we assume that 
\begin{subequations}\label{Def perturbations convergence}
\begin{equation}\label{Def gradient control perturbation}
\esssup \, \left\{\vert (- \gradth\zeta_{n}(\xthilde), 1)\vert:\xthilde\in \Gamma\right\} 
\leq C_{0}\quad \forall\, n\in \N.
\end{equation}
i.e. the gradients are globally bounded. Additionally assume uniform convergence 
\begin{equation}\label{Def uniform convergence perturbation}
 \Vert \zeta_{n} \Vert_{\scriptscriptstyle C(\Gamma)}  \xrightarrow[n\rightarrow \infty]{} 0
\end{equation}
\end{subequations}
%
%
From now on we denote $\Gamma^{n} = \Gamma^{\zeta_{n}}$ and $q^{n} = q^{\zeta_{n}}$.

For the forcing terms we assume there exists an open set $G$ containing $\{\Gamma^{n}\}$ and an element $\Phi \in \Hdiv (G)$ with $G$ an open region such that $\Phi\cdot\n\vert_{\scriptscriptstyle \Gamma^{n}} = f$ for all $n$. Here $\n$ denotes the upwards normal vector to $\Gamma^{n}$ and 
\begin{equation*}
\Hdiv= \{\v \in \Ltwo(G): \div \v \in L^{2}(G)\}.
\end{equation*}
Define 
\begin{equation}\label{Def trapped in between domain}
   U^{n} \defining \bigcup_{\xthilde \in \Gamma} (0, \zeta_{n}(\xthilde)) 
   \cup (\zeta_{n}(\xthilde), 0)
\end{equation}
then, $\partial U^{n} = \Gamma^{n} \cup \Gamma$. Due to condition \eqref{Def gradient control perturbation} the domain $U^{n}$ has Lipschitz boundary, then the classical duality relationship \cite{Tartar} holds i.e.
\begin{equation*}
\int_{\Gamma^{n}} f\, q^{n}\, dS - \int_{\Gamma} f\,q^{n}\,dS
=\int_{\partial U^{n}} \Phi\cdot \boldsymbol{\widehat{\nu}} \; q^{n}\, dS\\
= \int_{U^{n}} \div \Phi \, q^{n} + \Phi\cdot \grad q^{n}
\end{equation*}
then
\begin{multline}\label{Ineq estimate of surface term on the diagonal}
\int_{\Gamma^{n}} f\, q^{n}\, dS \leq  \Vert f \Vert_{-\frac{1}{2}, \Gamma}\,
\Vert q^{n}\Vert_{\frac{1}{2}, \Gamma}
+ \Vert \Phi\Vert_{\scriptscriptstyle \Hdiv(\Gamma)} \Vert q^{n}\Vert_{1, \Omega}\\
\leq \{ \Vert f \Vert_{-\frac{1}{2}, \Gamma}
+ \Vert \Phi\Vert_{\scriptscriptstyle \Hdiv(\Gamma)} \} \Vert q^{n}\Vert_{1, \Omega}.
\end{multline}
Combining \eqref{Ineq estimate of surface term on the diagonal} with \eqref{Ineq first a-priori estimate} gives
\begin{equation}\label{Ineq global a-priori estimate}
\Vert \,q^{n}\Vert_{1,\,\Omega}
\leq \{ C_{\scriptscriptstyle \Omega}\Vert  F\Vert_{0, \Omega}  
+ \Vert f \Vert_{-\frac{1}{2}, \Gamma}
+ \Vert \Phi\Vert_{\scriptscriptstyle \Hdiv(\Gamma)} \} .
\end{equation}
Due to the Rellich-Kondrachov theorem, there must exist a subsequence, denoted $\{q^{k}\}$ and an element $q^{*}\in H^{1}(\Omega)$ such that  
\begin{equation*}
   q^{k}\rightarrow q^{*}\,\;\text{weekly in}\;H^{1}(\Omega) ,
   \,\;\text{strongly in}\;L^{2}(\Omega).
\end{equation*}
Denoting $\Omega_{i}^{k} = \Omega_{i}^{\zeta_{k} }$ for $i = 1,2$, the variational statement can be written as
\begin{equation}\label{Stmt perturbed problem on cv subsequence}
\int_{\Omega} \grad q^{k} \cdot\grad r \,\ind_{\Omega_{1}^{k} }
 + \frac{1}{\epsilon} \int_{\Omega} \grad q^{k}\cdot\grad r \,\ind_{\Omega_{2}^{k} }\\
= \int_{\Omega} F\, r + \int_{\Gamma^{k} } f\, r\, d\,S
\end{equation}
for $r\in V$ arbitrary. Lett $\zeta_{k}\rightarrow 0$ in $C(\Gamma)$, first observe that for any $r\in H^{1}(\Omega)$ holds
\begin{equation*}
\left \vert \int_{\Gamma^{k}} f\, r\, dS - \int_{\Gamma} f\,r\,dS \right \vert
%
= \left \vert \int_{U^{k}} \div \Phi \, r + \Phi\cdot \grad r \right \vert
\leq \Vert \Phi \Vert_{\scriptscriptstyle \Hdiv(U^{k})} \Vert r \Vert_{\scriptscriptstyle H^{1}(U^{k})}
\end{equation*}
Since the right hand side converges to $0$ as $k\rightarrow \infty$ it follows that $ \int_{\Gamma^{k}} f\, r\, dS \rightarrow \int_{\Gamma} f\,r\,dS$. 
Next observe that
$\{\grad r\,\ind_{\Omega_{i}^{k}} \}$ converges strongly in $\textbf{L}^{\!2}(\Omega)$; together with the convergence of the surface forcing terms previously discussed, the expression \eqref{Stmt perturbed problem on cv subsequence} converges to
\begin{equation*}
\int_{\Omega} \grad q^{*} \cdot\grad r \,\ind_{\Omega_{1} }
 + \frac{1}{\epsilon} \int_{\Omega} \grad q^{*}\cdot\grad r \,\ind_{\Omega_{2}}\\
= \int_{\Omega} F\, r + \int_{\Gamma} f\, r\,dS.
\end{equation*}
Since $q^{*}$ is in $V$ and the variational statement above holds for all $r\in V$ the uniqueness of problem \eqref{Pblm Direct Formulation Original} implies $q^{*} = p$. The reasoning above holds for any subsequence of $\{q^{n}\}$ and the solution of \eqref{Pblm Direct Formulation Original}  is unique, then it follows that the whole sequence converges to $p$ i.e.
\begin{equation}\label{Stmt weak convergence of the pressures}
   q^{n}\rightarrow p\,\;\text{weekly in}\;H^{1}(\Omega) ,
   \,\;\text{strongly in}\;L^{2}(\Omega) .
\end{equation}
We close the section with an important observation. Test the statements \eqref{Pblm Direct Formulation Perturbed} on the diagonal $q^{n}$ and let $n\rightarrow \infty$; it yields
\begin{multline}\label{Eq best possible sense of norm convergence}
\lim_{n \rightarrow \infty} \left\{\int_{\Omega_{1}^{n}} \vert  \grad q^{n} \vert^{2}
 + \frac{1}{\epsilon} \int_{\Omega_{2}^{n}} \vert \grad q^{n}\vert^{2} \right\} 
= \int_{\Omega} F\, p + \int_{\Gamma} f\, p\,dS \\
= \int_{\Omega_{1}} \vert  \grad p \vert^{2}
 + \frac{1}{\epsilon} \int_{\Omega_{2} } \vert \grad p\vert^{2}.
\end{multline}
The map $r\mapsto \{\int_{\Omega_{1}} \vert  \grad r \vert^{2}
 + \frac{1}{\epsilon} \int_{\Omega_{2} } \vert \grad r\vert^{2} \}^{\frac{1}{2}}$ is a norm equivalent to the norm $\Vert \cdot\Vert_{V}$. However, due to the presence of the domains $\Omega^{n}_{i}$, $i = 1, 2$ the equality \eqref{Eq best possible sense of norm convergence} is not a statement of norms convergence which, together with the weak convergence, would allow to conclude strong convergence. However, due to the weak convergence and the equivalence of the norms $r\mapsto \{\int_{\Omega_{1}} \vert  \grad r \vert^{2}
  + \frac{1}{\epsilon} \int_{\Omega_{2} } \vert \grad r\vert^{2} \}^{\frac{1}{2}}$ and $\Vert \cdot \Vert_{V}$ we can conclude that
\begin{equation} \label{Ineq lim inf inequality}
 \int_{\Omega_{1}} \vert  \grad p \vert^{2}
 + \frac{1}{\epsilon} \int_{\Omega_{2} } \vert \grad p\vert^{2}
 \leq \liminf_{n} \left\{\int_{\Omega_{1}} \vert  \grad q^{n} \vert^{2}
  + \frac{1}{\epsilon} \int_{\Omega_{2}} \vert \grad q^{n}\vert^{2} \right\} .
\end{equation}
%
%
 
 %
 %
 %
 %
 %
 %
 %
 \subsection{The Strong Convergence}\label{Sec Strong Convergence}
Given a function $r\in V$ consider the following identities
\begin{subequations}\label{Eq Integration Strategy}
\begin{equation}\label{Eq Integration Strategy Domain 1}
\int_{\sOmrz}\vert \grad r \vert^{2} = \int_{\sOmrz - \Omega_{1}}\vert \grad r \vert^{2}  - \int_{\sOmfz - \Omega_{2}}\vert \grad r \vert^{2}  + \int_{\Omega_{1}}\vert \grad r \vert^{2} 
\end{equation}
\begin{equation}\label{Eq Integration Strategy Domain 2}
\int_{\sOmfz}\vert \grad r \vert^{2} = \int_{\sOmfz - \Omega_{2}}\vert \grad r \vert^{2}  - \int_{\sOmrz - \Omega_{1}}\vert \grad r \vert^{2}  + \int_{\Omega_{2}}\vert \grad r \vert^{2} 
\end{equation}
\end{subequations}
We define the perturbation term as
\begin{equation}\label{Eq perturbation term}
\varXi_{\zeta}(r) \defining \int_{\sOmfz - \Omega_{2}}\vert \grad r \vert^{2}  - \int_{\sOmrz - \Omega_{1}}\vert \grad r \vert^{2}  
\end{equation}
Moreover, the perturbation term satisfies that
\begin{multline}\label{Eq perturbation piece estimate}
\left\vert \varXi_{\zeta}(r) \right\vert= \left\vert\int_{\sOmfz - \Omega_{2}}\vert \grad r \vert^{2}  - \int_{\sOmrz - \Omega_{1}}\vert \grad r \vert^{2}  \right\vert \leq \left\vert\int_{\sOmfz - \Omega_{2}}\vert \grad r \vert^{2}   \right\vert + \left\vert \int_{\sOmrz - \Omega_{1}}\vert \grad r \vert^{2}  \right\vert \\
\leq \Vert r \Vert^{2}_{V}\left[\lambda_{N}(\sOmrz - \Omega_{1}) + \lambda_{N}(\sOmfz - \Omega_{2})\right]. 
\end{multline}
Therefore the following estimate holds
\begin{multline*} 
\int_{\sOmrz}\vert \grad r \vert^{2} + \frac{1}{\epsilon}\int_{\sOmfz}\vert \grad r \vert^{2} 
%
= \int_{\Omega_{1}}\vert \grad r \vert^{2}  + \frac{1}{\epsilon}\int_{\Omega_{1}}\vert \grad r \vert^{2}  + \left(1 - \frac{1}{\epsilon}\right)\varXi_{\zeta}(r) \\
\geq
\int_{\Omega_{1}}\vert \grad r \vert^{2}  + \frac{1}{\epsilon}\int_{\Omega_{1}}\vert \grad r \vert^{2}
- \left\vert 1 - \frac{1}{\epsilon}\right\vert \left[\lambda_{N}(\sOmrz - \Omega_{1}) + \lambda_{N}(\sOmfz - \Omega_{2})\right]  \Vert r\Vert_{V}^{2}.
\end{multline*}
We know $\Vert r\Vert_{V}^{2} \geq \epsilon \, \{\int_{\Omega_{1}} \vert  \grad r \vert^{2}
 + \frac{1}{\epsilon} \int_{\Omega_{2} } \vert \grad r\vert^{2}\}$ then, combining with the expression above we get
\begin{multline*} 
\int_{\sOmrz}\vert \grad r \vert^{2} + \frac{1}{\epsilon}\int_{\sOmfz}\vert \grad r \vert^{2} 
%
\\
\geq
\left(1 - \epsilon \left\vert 1 - \frac{1}{\epsilon}\right\vert \left[\lambda_{N}(\sOmrz - \Omega_{1}) + \lambda_{N}(\sOmfz - \Omega_{2})\right] \right) 
\left\{\int_{\Omega_{1}}\vert \grad r \vert^{2}  + \frac{1}{\epsilon}\int_{\Omega_{1}}\vert \grad r \vert^{2}\right\}.
\end{multline*}
Defining 
\begin{equation}\label{Eq lower bounding constant}
C_{\zeta} \defining
1 - \epsilon \left\vert 1 - \frac{1}{\epsilon}\right\vert \left[\lambda_{N}(\sOmrz - \Omega_{1}) + \lambda_{N}(\sOmfz - \Omega_{2})\right] ,
\end{equation}
We get the estimate
\begin{equation}\label{Ineq Inner Product Estimate}
\int_{\sOmrz}\vert \grad r \vert^{2} + \frac{1}{\epsilon}\int_{\sOmfz}\vert \grad r \vert^{2} 
%
%
\geq
(1 - C_{\zeta}) \left\{\int_{\Omega_{1}}\vert \grad r \vert^{2}  + \frac{1}{\epsilon}\int_{\Omega_{1}}\vert \grad r \vert^{2}\right\}
\quad \forall \, r\in V.
\end{equation}
In particular for the sequence of solutions $\{q^{n}: n\in \N\}\subseteq V$ holds
\begin{equation*} 
 C_{\zeta_{n}}\left\{\int_{\Omega_{1}} \vert  \grad q^{n} \vert^{2}
 + \frac{1}{\epsilon} \int_{\Omega_{2}} \vert \grad q^{n}\vert^{2} \right\}\leq 
\int_{\Omega_{1}^{n}} \vert  \grad q^{n} \vert^{2}
 + \frac{1}{\epsilon} \int_{\Omega_{2}^{n}} \vert \grad q^{n}\vert^{2} 
\end{equation*}
Letting $n\rightarrow \infty$ gives $\Vert \zeta_{n} \Vert_{C(\Gamma)}\rightarrow 0$ and consequently $C_{\zeta_{n}}\rightarrow0$. Then, taking $\limsup_{n}$ in the expression above yields
\begin{multline} \label{Ineq lim sup inequality}
\limsup_{n}\left\{\int_{\Omega_{1}} \vert  \grad q^{n} \vert^{2}
 + \frac{1}{\epsilon} \int_{\Omega_{2}} \vert \grad q^{n}\vert^{2} \right\}\\
 \leq 
\limsup_{n}\left\{\int_{\Omega_{1}^{n}} \vert  \grad q^{n} \vert^{2}
 + \frac{1}{\epsilon} \int_{\Omega_{2}^{n}} \vert \grad q^{n}\vert^{2} \right\} 
 = \int_{\Omega_{1}} \vert  \grad p \vert^{2}
  + \frac{1}{\epsilon} \int_{\Omega_{2} } \vert \grad p\vert^{2}.  
\end{multline}
Where the last equality holds due to \eqref{Eq best possible sense of norm convergence}. Putting together \eqref{Ineq lim inf inequality} and \eqref{Ineq lim sup inequality} we conclude 
\begin{equation} \label{Eq convergence of norms}
\lim_{n}\left\{\int_{\Omega_{1}} \vert  \grad q^{n} \vert^{2}
 + \frac{1}{\epsilon} \int_{\Omega_{2}} \vert \grad q^{n}\vert^{2} \right\}
 = \int_{\Omega_{1}} \vert  \grad p \vert^{2}
  + \frac{1}{\epsilon} \int_{\Omega_{2} } \vert \grad p\vert^{2}.  
\end{equation}
i.e. the norms $r\mapsto \{\int_{\Omega_{1}} \vert  \grad r \vert^{2}
  + \frac{1}{\epsilon} \int_{\Omega_{2} } \vert \grad r\vert^{2} \}^{\frac{1}{2}}$ converge and, due to the equivalence with the $V$-norm it follows that $\Vert q^{n}\Vert_{V} \vert \rightarrow \Vert p \Vert_{V}$. Since $q^{n}$ converges weakly to $p$ in $V$ it follows that 
\begin{equation} \label{Eq strong convergence}
\left\Vert   q^{n}  - p \right\Vert^{2}_{V}
\rightarrow 0\, \quad \text{as}\;\; n\rightarrow \infty. 
\end{equation}
%
%
%
%
%
%
\section{The One Dimensional Case}\label{Sec One dimensional example}
Here, we restrict our attention to the one dimensional problem in order to gain deep insight on the phenomenon. An example of how valuable this approach is, can be found in \cite{Gunzburger2009}. For the problem in one dimensional setting we choose $\Omega_{\,1} \defining (-1, 0)$, $\Omega_{\,2} \defining (0, 1)$ and the interface $\Gamma = \{0\}$. In this context a perturbation is given by a single point $\zeta$; the perturbed domains are given by $\Omrz \defining (-1, \zeta)$, $\Omfz\defining(\zeta, 1)$ and the perturbed interface $\Gammaz \defining \{\zeta\}$. Clearly $\Omega  = \Omega_{1}\cup\Gamma\cup\Omega_{2} = \Omrz\cup \Gammaz\cup \Omfz = (-1,1)$, see figure \ref{Fig Orthogonal Decomposition}. Notice that in the one dimensional case the space $V$ and its inner product given in definition \ref{Def function space setting} reduce to
\begin{subequations}\label{Def reduction ot 1-D function space}
   \begin{equation}\label{Def test space direct formulation}
V =\left\{r\in H^{1}(-1, 1): r(-1) = 0\right\} ,
\end{equation}
\begin{equation}\label{Def the right inner product}
\langle \pi,\, \kappa \rangle_{\scriptscriptstyle V} = \int_{-1}^{1}\del \pi\, \del \kappa .
\end{equation}
\end{subequations}
Where $\del$ indicates de weak derivative. Similarly, the problems \eqref{Pblm Direct Formulation Perturbed} and \eqref{Pblm Direct Formulation Perturbed} transform in 
\begin{subequations}\label{Eq direct 1-d}
\begin{equation}\label{Eq direct original 1-d}
p\in V:\quad
\int_{-1}^{\,0}\del p\;\del r
+ \frac{1}{\epsilon} \int_{0}^{\,1}\del p\;\del r\\
= \int_{-1}^{\,1}F\,r
+ f \left(0\right) r \left(0\right)\quad\forall\;r\in V.
\end{equation}
\begin{equation}\label{Eq direct perturbed 1-d}
q\in V:\quad
\int_{-1}^{\,\zeta}\del q\;\del r
+ \frac{1}{\epsilon} \int_{\zeta}^{1}\del q\;\del r\\
= \int_{-1}^{1}F\,r
+ f \left(\zeta\right) r \left(\zeta\right)\quad\forall\;r\in V.
\end{equation}
\end{subequations}
%
%
%
%
%
%
%
%
%
\subsection{The Subspace $H$ and its Orthogonal Projection}
\begin{figure}[!]
\caption[1]{Orthogonal Decomposition}\label{Fig Orthogonal Decomposition}\vspace{5pt}
\centerline{\resizebox{8cm}{8cm}
{\includegraphics{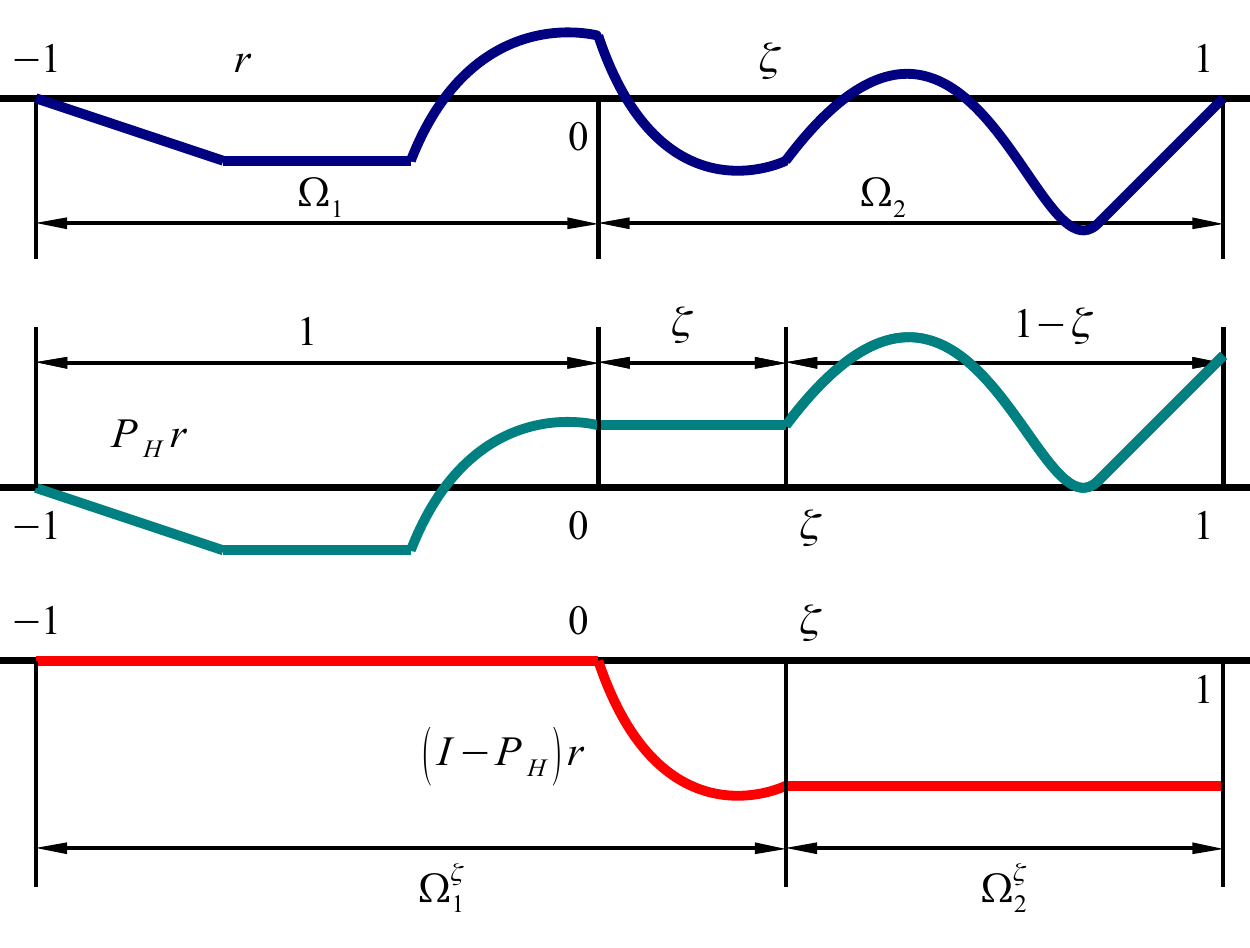}}}
\end{figure}

In order to estimate the norms $\Vert p - q^{\zeta}\Vert_{\,1,\Omega}$, $\Vert p - q^{\zeta}\Vert_{\,0,\Omega}$ we need to project the solutions $p$ and $q^{\zeta}$ into the adequate subspace using the convenient geometry defined by the inner product \eqref{Def the right inner product}. For simplicity, from now on it will be assumed $\zeta>0$. Consider the subspaces
\begin{subequations}\label{Def direct sum decomposition}
\begin{equation}\label{Def the right space}
H\defining\left\{\kappa\in V: \del \kappa = 0 \;\text{on}\;(0,\zeta)\right\}
\end{equation}
\begin{equation}\label{Def orthogonal the right space}
H^{\perp}\defining\left\{\pi\in V: \langle\pi, \kappa\rangle_{\scriptscriptstyle V} = 0
\;\forall\;\kappa\in H\right\}
\end{equation}
\end{subequations}
Next, we characterize the structure of $H^{\perp}$.
\begin{lemma}\label{Th characterization of Kperp}
Let $H^{\perp}$ and $H$ defined in \eqref{Def direct sum decomposition} then
\begin{equation}\label{Def orthogonal to the right space}
H^{\perp} = \left\{\pi\in V: \pi = 0\;\text{on}\;(-1,0)\right.\\
\left.\,,
\del \pi = 0\;\text{on}\;(\zeta, 1)\right\}
\end{equation}
\begin{proof} It is direct to see that if $\pi\in V$ is such that $\pi = 0$ in $(-1,0)$ and $\del\pi = 0$ in $(\zeta,1)$ then $\pi\in H^{\perp}$. For the other inclusion take $\rho\in C_{\,0}^{\,\infty}\left(-1, \,0\right)$ such that $\int_{-1}^{\,\0}\rho\,dx = 1$ and extended it by zero to the whole domain $\left(-1, \,1\right)$. Choose any $\phi\in C_{\,0}^{\infty}\left(-1, \,0\right)$, extended it by zero to $\left(-1,\,1\right)$ and build the auxiliary function
\begin{equation*}
\Phi(x) \defining 
\int_{-1}^{\,x}\phi(t)\,d\,t\,\ind_{(-1, 0)}\, (x)\\
- \int_{-1}^{\,0}\phi(y)\,d y \int_{-1}^{\,x} \rho (t)\,d\,t\,\ind_{(-1, 0)}\, (x) .
\end{equation*}
It is direct to see that $\Phi\in H$. Now take any $\pi\in H^{\perp}$, then
\begin{multline*}
0 = \left\langle \pi,\,\Phi\right\rangle_{\scriptscriptstyle V} = \int_{-1}^{\,0}\del \pi\; \del \Phi 
=
\int_{-1}^{\,0}\del \pi\,(x)\,\left[\phi\,(x)-\left(\int_{-1}^{\,0}\phi(y)\,d y\right)\rho(x)\right] dx\\
= \int_{-1}^{\,0}\del \pi\,(x)\,\phi\,(x) \,dx
-\int_{-1}^{\,0}\phi(y) \,dy
\int_{-1}^{\,0}\del \pi\,(x)\,\rho(x)\, dx .
\end{multline*}
i.e.
\begin{equation*}
\int_{-1}^{\,0}\del \pi\,(x)\,\phi\,(x)\, dx
=\int_{-1}^{\,0}\phi(y)\, dy\int_{-1}^{\,0}\del \pi\,(x)\,\rho(x)\,dx ,
\end{equation*}
for all $\phi\in C_{\,0}^{\,\infty}\left(-1,\,0\right)$. Therefore, we conclude $\del  \pi$ must be constant in $\left(-1,\,0\right)$. Using an analogous construction we also conclude $\del \pi$ must be constant in $\left(\zeta,\,1\right)$. Now we prove that such constants must be zero. Consider any $\kappa \in H$, then we have
\begin{multline*}
0 = 
\left\langle \pi,\,\kappa\right\rangle_{\scriptscriptstyle V} 
=  \int_{-1}^{\,0}\del \pi\,\del \kappa
+\int_{\zeta}^{\,1}\del \pi\,\del \kappa 
=\del \pi(-1)\int_{-1}^{\,0}\del \kappa
+
\del \pi(1)\int_{\zeta}^{\,1}\del \kappa
\\
= \del \pi(-1)\,\left(\kappa\left(0\right) - \kappa\left(-1\right)\right)
+ \del \pi(1)\,\left(\kappa\left(1\right) - \kappa\left(\zeta\right)\right) .
\end{multline*}
Since $\kappa\in H\subset V$ it holds $\kappa\left(-1\right) = 0$ and the above expression writes
\begin{equation}\label{Eq characterization of del u}
\del \, \pi(-1)\,\kappa\left(0\right)
+ \del \, \pi(1)\,\left(\kappa\left(1\right) - \kappa\left(\zeta\right)\right) = 0 .
\end{equation}
Due to $\del \kappa = 0 $ on $\left(0, \,\zeta\right)$ the function $\kappa$ must be constant on this interval, therefore $\kappa\left(0\right) = \kappa\left(\zeta\right)$. Recalling the above holds for any $\kappa\in H$, choose a test function such that $\kappa\left(0\right) = \kappa\left(\zeta\right) = 0$ and $\kappa(1) \neq 0$, then \eqref{Eq characterization of del u} reduces to $\del \pi(1)\,\kappa\left(1\right) = 0$ and we conclude $\del \pi (1)  = 0$. Hence, \eqref{Eq characterization of del u} reduces to $\del \pi (-1)\,\kappa (0) = 0$. Since $\kappa\in H$ is arbitrary we know $\kappa(0)$ need not be zero for all $\kappa\in H$, then we conclude $\del \pi (-1) = 0$. Therefore $\pi$ must be constant on the intervals $\left(-1, \,0\right)$ and $\left(\zeta,\, 1\right)$. Finally, the fact that $\pi\in H^{\,\perp}\subset V$ yields $\pi(-1) = 0$; this implies $\pi = 0$ on $\left(-1,\,0\right)$ which completes the proof.
\end{proof}
\end{lemma}
Now we present the characterization of the orthogonal projections onto the subspaces $H$ and $H^{\perp}$.
\begin{theorem}\label{Th characterization of projections}
Let $H, H^{\perp}$ be the spaces defined in \eqref{Def direct sum decomposition}. Denote $P_{\scriptscriptstyle H}$ and $P_{\scriptscriptstyle H^{\perp}}$ the orthogonal projections onto the subspaces $H$ and $H^{\perp}$ respectively. Then, for any $r\in V$ holds
\begin{subequations}\label{Eq orthogonal decomposition}
\begin{equation}\label{Eq orthogonal projection on W}
P_{\scriptscriptstyle H}\,r (x)
= r (x)  \ind_{\left[-1,\,0\,\right]}(x)
+ r\left(0\right) \ind_{\left[\,0,\,\zeta\,\right]}(x) 
+\left\{ r(x) - \left[ r\left(\zeta\right) - r\left(0\right)\right]\right\} \ind_{\left[\zeta,\,1\right]}(x)
\end{equation}
\begin{equation}\label{Eq orthogonal projection on W perp}
P_{\scriptscriptstyle  H^{\perp}} r (x) =  
\left [r(x) - r\left(0\right)\right] \ind_{\left[\,0,\,\zeta\,\right]}(x)\\
+ \left[ r\left(\zeta\right) - r\left(0\right)\right] \ind_{\left[\zeta,\,1\right]}(x)
\end{equation}
\end{subequations}
Where $ \ind_{A}(\cdot)$ denotes the indicator function of the set $A$.

\begin{proof} For any $r\in V$ it is direct to see that the function $x\mapsto r (x)  \ind_{\left[-1,\,0\,\right]}(x) + r\left(0\right) \ind_{\left[\,0,\,\zeta\,\right]}(x) +\left\{ r(x) - \left[ r\left(\zeta\right) - r\left(0\right)\right]\right\} \ind_{\left[\zeta,\,1\right]}(x)$ is in $H$ and that the map $x\mapsto \left [r(x) - r\left(0\right)\right]\,\ind_{\left[\,0,\,\zeta\,\right]}(x)
+ \left[ r\left(\zeta\right) - r\left(0\right)\right] \,\ind_{\left[\zeta,\,1\right]}(x)$ belongs to $H^{\perp}$. Also, their sum gives $r$. The result follows due to the characterization given in lemma \ref{Th characterization of Kperp}.
\end{proof}
\end{theorem}
\begin{remark}
In order to better understand the nature of the orthogonal decomposition we present figure \ref{Fig Orthogonal Decomposition}. An absolutely continuous function $r$ (blue line) is decomposed in $P_{\scriptscriptstyle H}\,r$ (turquoise line) and $P_{\scriptscriptstyle H^{\perp}} r =  (I - P_{\scriptscriptstyle H} ) r$ (red line). Also, the domains $\Omega_{1}, \Omega_{2}$ and $\Omega_{1}^{\zeta}, \Omega_{2}^{\zeta}$ are depicted.
\end{remark}
%
%
%
%
%
%
%
%
%
%
\subsection{The Problems Restricted to $H$}
Test the problem \eqref{Eq direct original 1-d} with a function $\kappa\in H$, we have
\begin{equation*}
\int_{-1}^{\,0}\del p\, \del \kappa
+ \frac{1}{\epsilon}\int_{\zeta}^{\,1} \del p\, \del \kappa =
\int_{-1}^{\,1}F\, \kappa
+ f\left(0\right) \kappa\left(0\right) .
\end{equation*}
Now decompose $p$ in $P_{\scriptscriptstyle H}p$ and $P_{\scriptscriptstyle H}p$ using \eqref{Eq orthogonal projection on W} and \eqref{Eq orthogonal projection on W perp}; we get
\begin{equation*}
%
%
\int_{-1}^{\,0}\del\, \left(P_{\scriptscriptstyle H}\,p\right)\, \del \kappa
+ \frac{1}{\epsilon}\int_{\zeta}^{\,1} \del  \left(P_{\scriptscriptstyle H}\,p\right)\, \del \kappa 
%
%
= \int_{-1}^{\,0}F\, \kappa + \int_{\zeta}^{\,1}F\, \kappa 
+ \left[\int_{0}^{\,\zeta}F
+ f\left(0\right) \right]\kappa\left(0\right) .
\end{equation*}
Here, the last equality used the fact that $\kappa$ is constant in $\left(0,\, \zeta\right)$ for all $\kappa\in H$. We write the statement as
\begin{multline}\label{Eq ueps on W}
P_{\scriptscriptstyle H} p\in H: \int_{-1}^{\,0}\del\, \left(P_{\scriptscriptstyle H}\, p\right)\, \del \kappa 
+ \frac{1}{\epsilon}\int_{\zeta}^{\,1} \del  \left(P_{\scriptscriptstyle H}\, p\right)\, \del \kappa \\
=
\int_{-1}^{\,0}F\, \kappa
+ \int_{\zeta}^{\,1}F\, \kappa 
+ \left[\int_{0}^{\,\zeta}F
+ f\left(0\right) \right]\kappa\left(0\right) , \;\forall\, \kappa\in H .
\end{multline}
On the other hand consider the problem
\begin{multline}\label{Eq original problem on W}
\sigma\in H: \int_{-1}^{\,0}\del \sigma\; \del \kappa
+ \frac{1}{\epsilon}\int_{\zeta}^{\,1} \del  \sigma\; \del \kappa\\
=
\int_{-1}^{\,0}F\, \kappa
+ \int_{\zeta}^{\,1}F\, \kappa 
+ \left[\int_{0}^{\,\zeta}F
+f\left(0\right) \right]\kappa\left(0\right) , \; \forall\, \kappa\in H .
\end{multline}
The bilinear the form $\mathcal{A}\left(\pi, \kappa\right)\defining\int_{-1}^{\,0}\del\, \pi\, \del \kappa + \frac{1}{\epsilon}\int_{\,\zeta}^{\,1} \del  \pi\, \del \kappa$ is $H$-elliptic and continuous, therefore the problem \eqref{Eq original problem on W} is well-posed and we conclude that $P_{\scriptscriptstyle H} p$ is the unique solution of \eqref{Eq original problem on W}.
%
%
%
%
Repeating the same procedure on the perturbed problem \eqref{Eq direct perturbed 1-d} we conclude $P_{\scriptscriptstyle H}q$ is the unique solution to the well-posed variational problem
%
%
\begin{multline}\label{Eq mueps on W}
P_{\scriptscriptstyle H} q\in H: \int_{-1}^{\,0}\del\, \left(P_{\scriptscriptstyle H}\, q\right)\, \del \kappa 
+ \frac{1}{\epsilon}\int_{\zeta}^{\,1} \del  \left(P_{\scriptscriptstyle H}\, q\right)\, \del \kappa \\
=
\int_{-1}^{\,0}F\, \kappa\, d x + \int_{\zeta}^{\,1}F\, \kappa 
+ \left[\int_{0}^{\,\zeta}F 
+ f\left(\zeta\right) \right]\kappa\left(\zeta\right) ,\;\forall\, \kappa\in H .
\end{multline}
%
%
%
%
%
%
\subsection{The Problems Restricted to $H^{\perp}$}  
We repeat the same strategy of the previous section and get
\begin{equation}\label{Eq ueps on W perp}
P_{\scriptscriptstyle H^{\perp}} p\in H^{\perp}: \frac{1}{\epsilon}\int_{0}^{\,\zeta}\del\, \left(P_{\scriptscriptstyle H^{\perp}}\, p\right)\, \del \kappa \\
=
\int_{0}^{\,\zeta}F\, \kappa
+ \int_{\zeta}^{\,1}F \, \kappa\left(\zeta\right) , \; \forall\, \kappa\in H^{\perp} .
\end{equation}
%
%
This is the weak solution of following strong problem
%
%
%
%
%
%
\begin{equation}\label{Eq magnified strong problem on W perp}
\begin{split}
-\del\,\frac{1}{\epsilon}\,\del P_{\scriptscriptstyle H^{\perp}}p&=F  \text{ in } \left(0,\, \zeta\right) ,\\[3pt]
P_{\scriptscriptstyle H^{\perp}} p  & = 0 \text{ in } \left[-1,\,0\right) ,\\[3pt]
P_{\scriptscriptstyle H^{\perp}} p&= \mathrm{constant}\,\text{ in } \left(\zeta,\, 1\right) ,\\
\frac{1}{\epsilon}\,\del P_{\scriptscriptstyle H^{\perp}} p\left(\zeta^{-}\right)
&= \int_{\zeta}^{\,1}F .
%
\end{split}
\end{equation}
Where $\del P_{\scriptscriptstyle H^{\perp}} p\left(\zeta^{-}\right) = \lim_{\,t\,\rightarrow\, \zeta^{-}} \del P_{\scriptscriptstyle H^{\perp}} p\left(t\right) $. In the same fashion
\begin{equation}\label{Eq mueps on W perp}
P_{\scriptscriptstyle H^{\perp}} q\in H^{\perp}: 
\int_{0}^{\,\zeta}\del\, \left(P_{\scriptscriptstyle H^{\perp} } q\right)\, \del \kappa \\
=
\int_{\,0}^{\,\zeta}F\, \kappa +
\left[\int_{\zeta}^{\,1}F 
+ f\left(\zeta\right) \right]\kappa (\zeta) , \;
\forall\, \kappa\in H^{\perp} .
\end{equation}
Where the simplification on the term of the right hand side has been made since $\kappa \left(x\right) = \kappa \left(\zeta\right)$ for $x\in \left(\zeta,\,1\right)$. Thus $P_{\scriptscriptstyle H^{\perp}} q$ is the solution to the strong problem
%
%
%
%
%
%
\begin{equation}\label{Eq magnified perturbed strong problem on W perp}
\begin{split}
-\del\,\del P_{\scriptscriptstyle H^{\perp}} q&=F  \text{ in } \left(0,\, \zeta\right) ,\\[3pt]
P_{\scriptscriptstyle H^{\perp}} q &=0 \text{ in } \left[-1,\,0\right) ,\\[3pt]
P_{\scriptscriptstyle H^{\perp}} q& = \mathrm{constant}\,\text{ in } \left(\zeta,\, 1\right) ,\\
\del P_{\scriptscriptstyle H^{\perp}} q\left(\zeta^{-}\right)
&= \int_{\zeta}^{\,1}F + f\left(\zeta\right) .
%
\end{split}
\end{equation}
%
%
%
%
%
%
\subsection{ Estimates for the $H$ Projections}  
Test \eqref{Eq ueps on W} and \eqref{Eq mueps on W} with $P_{\scriptscriptstyle H} p - P_{\scriptscriptstyle H} q$ and subtract the result to get
\begin{multline*}
\int_{-1}^{\,0}\del (P_{\scriptscriptstyle H} p - P_{\scriptscriptstyle H}  q )\, 
\del (P_{\scriptscriptstyle H} p - P_{\scriptscriptstyle H } q ) 
+ \frac{1}{\epsilon}\int_{\zeta}^{\,1} \del (P_{\scriptscriptstyle H}  p
- P_{\scriptscriptstyle H}\, q)\, 
\del (P_{\scriptscriptstyle H}  p-P_{\scriptscriptstyle H} q)\\
=
f\left(0\right) \left[P_{\scriptscriptstyle H} p (0) - P_{\scriptscriptstyle H} q (0)\right] 
- f\left(\zeta\right)\left[P_{\scriptscriptstyle H}  p (\zeta) - P_{\scriptscriptstyle H} q(\zeta)\right] \\[5pt]
=\left[f\left(0\right)-f\left(\zeta\right)\right] \left[P_{\scriptscriptstyle H} p (0)-P_{\scriptscriptstyle H} q(0)\right].
\end{multline*}
The last equality holds true since $\kappa\left(0\right) = \kappa\left(\zeta\right)$ for all $\kappa\in H$. Since $\vert r(x)\vert \leq \sqrt{2}\,\Vert \del r \Vert_{0,(-1, 1)}$ for all $r\in V$. Thus, the expression above can be estimated by 
%
%
%
%
\begin{equation}\label{Ineq Estimate on the projection onto W}
\Vert\, P_{\scriptscriptstyle H} p- P_{\scriptscriptstyle H}\, q\Vert_{\scriptscriptstyle V}
\leq C\,\vert f \left(0\right) - f (\zeta) \vert
\end{equation}
%
%
%
%
%
%
\subsection{Estimates for the $H^{\perp}$ Projections}  
Since \eqref{Eq magnified strong problem on W perp} and \eqref{Eq magnified perturbed strong problem on W perp} are both ordinary differential equations, the exact solutions can be found. These are given by
\begin{multline}\label{Eq exact ueps on W perp}
P_{\scriptscriptstyle H^{\perp}} p (x) 
= -\epsilon \int_{0}^{\,x}\!\!\!\int_{0}^{\,t} F(s)\, ds\, dt \,\ind_{\left[0,\,\zeta\,\right]}(x)\\
+\epsilon\, x \! \int_{0}^{\,1}\!\!\! F\,\ind_{\left[0,\,\zeta\,\right]}(x) 
-\epsilon \, \zeta \left[\int_{0}^{\,\zeta}\!\!\!\int_{0}^{\,t} \!\!\! F(s)\, ds\, dt
- \int_{0}^{\,1}\!\!\!F\right]\!\!  \ind_{\left[\zeta,\,1\right]}(x)
\end{multline}
\begin{multline}\label{Eq exact mueps on W perp}
P_{\scriptscriptstyle H^{\perp}} q (x) = 
- \int_{0}^{\,x} \!\!\! \int_{0}^{\,t} F(s)\, ds\, dt \, \ind_{\left[0,\,\zeta\,\right]}(x) 
+ \left[\int_{0}^{\,1}F
+ f(\zeta) \right] x\, \ind_{\left[0,\,\zeta\,\right]}(x)\\
- \left\{\int_{0}^{\zeta} \!\!\! \int_{0}^{\,t} F(s)\, ds\, dt
-\left[\int_{0}^{\,1}F + f(\zeta) \right]
\zeta\,\right\}\ind_{\left[\zeta,\,1\right]}(x)
\end{multline}
We estimate the norm of the difference using the exact expressions \eqref{Eq exact ueps on W perp} and \eqref{Eq exact mueps on W perp}.
%
%
%
%
When computing the $L^{\,2}$-norm of the difference of derivatives we get 
\begin{multline*}
\left\Vert  P_{\scriptscriptstyle H^{\perp}} p - P_{\scriptscriptstyle H^{\perp}} q \right\Vert_{\,V} =
\Vert \del P_{\scriptscriptstyle H^{\perp}} p -\del P_{\scriptscriptstyle H^{\perp}} q \Vert_{\scriptscriptstyle L^{2}\left(0, \zeta\right)} 
\leq \left(1 - \epsilon\right)
\left\Vert\int_{0}^{\,(\cdot)} F(t)\,  dt\,\ind_{\left[0,\,\zeta\,\right]}(\cdot)\right\Vert_{\scriptscriptstyle L^{2}\left(0,\, \zeta\right)} \\
+\sqrt{\zeta} \left[ (1 - \epsilon) \int_{0}^{1}F(t)\, dt
+ f(\zeta)\right] 
%
%
\leq \frac{\zeta (1-\epsilon)}{\sqrt{2}}\,\Vert F\Vert_{\scriptscriptstyle L^{2}(0, \, \zeta)} 
+ \zeta \left\vert f(\zeta) + ( 1 - \epsilon ) \int_{0}^{1}F \right\vert .
\end{multline*}
Then we conclude
\begin{equation}\label{Ineq Estimate on the projection onto W perp}
\Vert  P_{\scriptscriptstyle H^{\perp}} p - P_{\scriptscriptstyle H^{\perp}} q \Vert_{\scriptscriptstyle V}
\leq C \,\zeta \,\{\Vert F \Vert_{ \scriptscriptstyle L^{2}(-1, \,1)} 
+ \vert f(\zeta)\vert \}
\end{equation}
Where $C>0$ is an adequate constant.
%
%
%
%
%
%
%
\subsection{Global Estimate of the Perturbation} \label{Sec global estimate of perturbation}
For the global estimate of the difference recall $\Vert \y \Vert_{2} \leq \Vert\y\Vert_{1}$ for each $\y\in \R^{2}$ and combine the estimates \eqref{Ineq Estimate on the projection onto W},  \eqref{Ineq Estimate on the projection onto W perp} to get
\begin{multline}\label{Ineq continuous dependence wrt perturbation}
\Vert\,p-q\Vert_{\scriptscriptstyle V} 
=\{ \Vert P_{\scriptscriptstyle H}( p - q )\Vert_{\scriptscriptstyle V}^{2}
+\Vert P_{\scriptscriptstyle H^{\perp}}(p - q)\Vert_{\scriptscriptstyle V}^{2} \}^{\frac{1}{2}} 
\leq \Vert\,P_{\scriptscriptstyle H} p 
- P_{\scriptscriptstyle H} q \Vert_{\scriptscriptstyle V}
+ \Vert\,P_{\scriptscriptstyle H^{\perp}} p 
- P_{\scriptscriptstyle H^{\perp}} q \Vert_{\scriptscriptstyle V} \\[4pt]
\leq \sqrt{2} \, \vert f(0) - f(\zeta) \vert
+ C\,\zeta \, \{ \Vert F\Vert_{\scriptscriptstyle L^{\,2} (-1, \,1) }
+ \vert f(\zeta)\vert \} .
\end{multline}
In order to have continuous dependence of the solutions with respect to perturbations of the interface, it is direct to see that the finiteness of $\left\Vert F\right\Vert_{\scriptscriptstyle L^{\,2} (-1, \,1)}$ needs to be required, and that conditions on the forcing term $f$ behavior need to be stated. In the last line of inequality \eqref{Ineq continuous dependence wrt perturbation} the third summand needs $f$ to be bounded in a neighborhood $[0, \delta)$ while the first summand demands it to be right continuous in a neighborhood $[0, \delta)$ for some $\delta >0$. Recalling $H_{\text{div}}(0, \delta) = H^{1}(\delta)$ in one dimension, the hypothesis that $f = \del \Phi$ for some $\Phi\in H_{\text{div}}(0, \delta)$ assumed in section \ref{Sec stability estimates} is sufficient to satisfy these conditions, however, it is not necessary.

Finally, a repetition of the the same procedure for perturbations to the left, i.e. when $\zeta <0$ yields
\begin{equation}\label{Ineq continuous dependence wrt left perturbation}
\left\Vert\,p - q\right\Vert_{\scriptscriptstyle V}
\leq \sqrt{2}\,\vert f(0) - f(\zeta) \vert \\
+ C \,\vert\zeta \vert \{ \left\Vert F\right\Vert_{\scriptscriptstyle L^{\,2}(-1, \,1)}
+ \left\vert f(0)\right\vert \}
\end{equation}
The section summarizes in the following result
\begin{theorem}\label{Th strong convergence in the one dimensional case}
   Let $F\in L^{2}(-1,1)$ and $f\in C(-\delta, \delta)$ for some $\delta>0$, then
   \begin{equation}
      \Vert p - q^{\zeta}\Vert_{\scriptscriptstyle V} \xrightarrow[\zeta \rightarrow 0]{} 0.
   \end{equation}
   \emph{i.e.} the sequence of perturbed solutions $\{q^{\zeta}\}$ converge strongly to the original one.
   \begin{proof}
      The estimate \eqref{Ineq continuous dependence wrt left perturbation} together with the hypotheses on the forcing terms yield the desired results.
   \end{proof}
\end{theorem}
\begin{remark}\label{Rem 1-D analysis projections}
   It is important to stress that the successful analysis in the one dimensional case, heavily relies on the dimension itself. The characterization of the right space $H$ and its orthogonal projections can not be done in a multiple dimensional setting, even in a very simple geometric domain such as the unit ball or the unit square. Also, the solutions provided by equations \eqref{Eq exact ueps on W perp}, \eqref{Eq exact mueps on W perp} are possible only due to one dimensional framework.
\end{remark}
\begin{remark}\label{Rem 1-D analysis dependence}
   The estimates \eqref{Ineq continuous dependence wrt perturbation} and \eqref{Ineq continuous dependence wrt left perturbation} heavily depend on the pointwise behavior of $f$, e.g. if $f(x) \equiv C \, x^{s}$ for $s > 0$ very small, the rate of convergence is very slow. It also reveals the nonlinear dependence of the solution $q$ with respect to $\zeta$.
\end{remark}
%
%
%
%
%
%
%
%
%
%
%
%
%
%
%
\section{A Simple Geometry in Multiple Dimensional Setting}\label{Sec comparing problems}
\begin{figure}[!]
\caption[1]{Flattening of the Interface}\label{Fig Flattening Interface}
\centerline{\resizebox{8cm}{7cm}
{\includegraphics{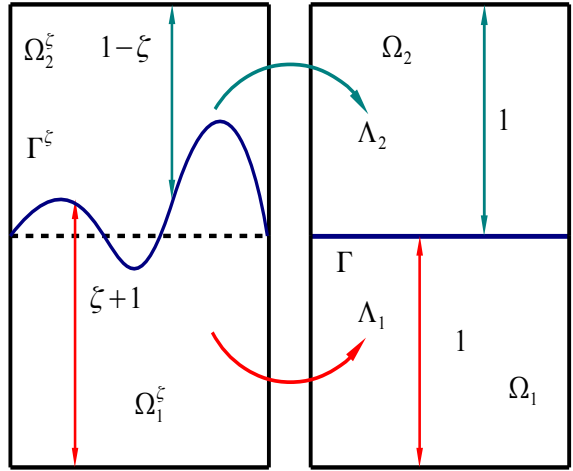}}}
\end{figure}

In this section we choose the simplest possible geometry in multiple dimensions in order to illustrate the nonlinearities that the phenomenon of interface geometric perturbation involves. Let $\Gamma\subseteq \R^{\! N-1}$ be open connected and $\Omega_{1} = \Gamma\times (0,1)$ and $\Omega_{2} = \Gamma\times(-1, 0)$ i.e. the domain on the right of figure \ref{Fig Flattening Interface}. Also assume that the perturbation $\zeta$ is piecewise $C^{1}(\Gamma)$.
We exploit the geometry defining the fractional bijective maps $\Lambda_{i}:\Omega_{i}^{\zeta}\rightarrow \Omega_{i}$ for $i = 1, 2$ as follows
\begin{subequations}\label{Def conformal map}
\begin{equation}\label{Def conformal map domain i}
\Lambda_{i}(\Xthilde, X_{\scriptscriptstyle N}) \defining (\Xthilde, \frac{X_{\scriptscriptstyle N} - \zeta(\Xthilde)}{1 - (-1)^{i} \zeta(\Xthilde)}) 
\end{equation}
Also define
\begin{equation}\label{Def global conformal map}
\Lambda(\X) \defining \sum_{i = 1, 2}
\Lambda_{i}(\X)  \ind_{\Omega_{i}}(\X) 
\end{equation}
\end{subequations}
%
%
%
Now set the new variables
\begin{align}\label{Def change of variables}
z\defining \Lambda(\Xthilde, X_{\scriptscriptstyle N})\cdot\enversor ,\\
\xthilde \defining \Lambda(\Xthilde, X_{\scriptscriptstyle N}) - \left(\Lambda(\Xthilde, X_{\scriptscriptstyle N})\cdot\enversor\right) \enversor .
\end{align}
%
%
%
%
%
%
\subsection{Changes on Gradient Structure}
For the maps $\Lambda_{i}$ denote $\Lambda_{i}'$ its derivative or Jacobian matrix, then for $i = 1, 2$ holds
\begin{subequations}\label{Eq derivatives domain 1}
\begin{equation}\label{Eq Jacobian domain 1}
\Lambda_{i}'
= \left[\begin{array}{cc}
I & \0\\[5pt]
(1 - (-1)^{i} z)\gradth^{T} \zeta  & 1 - (-1)^{i} \zeta
\end{array}
\right] .
\end{equation}
Since $\Vert \zeta \Vert_{\scriptscriptstyle C(\Gamma)} < 1$ for functions belonging to $\pert(\Gamma, \Omega)$ (defined in equation \eqref{Def perturbations of interface}), the absolute value of the determinant of the Jacobian matrix is given by
\begin{equation}\label{Eq det Jacobian domain 1}
\left\vert \det \left(\Lambda_{i}'\right)
\right\vert
= \left\vert 1 - (-1)^{i}\zeta(\xthilde\,)\right\vert = 1 - (-1)^{i}\zeta(\xthilde\,).
\end{equation}
\end{subequations}
For a scalar function we observe that whenever $X\in \Omega_{i}^{\zeta}$ the gradient has the following structure
\begin{subequations}\label{Eq gradient structure domain 1}
\begin{equation}\label{Eq extended gradient structure domain 1}
\left\{\begin{array}{c}
\gradth_{\!\! \scriptscriptstyle X}\\[5pt]
\dfrac{\partial}{\partial X_{\scriptscriptstyle N}}
\end{array}\right\} 
%
%
= \left[\begin{array}{cc}
I & -(-1)^{i}\,\dfrac{z - (-1)^{i} }{1- (-1)^{i} \zeta}\,\gradth \zeta\\[10pt]
\0^{T} & \dfrac{1}{1 - (-1)^{i} \zeta}
\end{array}
\right]
\left\{\begin{array}{c}
\gradth_{\!\! x}\\[10pt]
\dfrac{\partial}{\partial z}
\end{array}\right\} ,
\end{equation}
in matrix notation
\begin{equation}\label{Eq matrix notation gradient structure domain 1}
\grad_{\!\! \scriptscriptstyle X} = A_{i}^{\zeta}\,\grad_{\!\! x} \quad \text{for all} \, X\in \Omega_{i}^{\zeta}\,\; \text{and} \,\; i = 1, 2.
\end{equation}
\end{subequations}
%
%
%
%
%
%
%
%
\subsection{Fractional Mapping and the $H^{1}(\Omega)$ Space}\label{Sec Conformal maps and H1 spaces}
From now on we endow the set $\pert(\Omega, \Gamma)$ with the norm $W^{1, \infty}(\Gamma)$ i.e. the sum of the essential suprema for the function and its gradient. Define the following change of variable
\begin{definition}\label{Def morphism of functions}
For each element $r\in H^{1}(\Omega)$ we define the \emph{Fractional Mapping} operator by
\begin{equation}\label{Eq morphism of functions}
T\,r \defining \sum_{i = 1, 2}
\left(r\circ\Lambda_{i}^{-1}\right) \ind_{\Omega_{i}} 
\end{equation}
\end{definition}
%
%
%
%
\begin{lemma}\label{Th weak derivatives change of variable}
Let $r\in H^{1}(\Omega)$ then for each $1\leq \ell\leq N$ holds
%
%
\begin{equation}\label{Eq weak derivatives}
\frac{\partial}{\partial \,x_{\ell}}\,T r  = \sum_{i = 1, 2}
\frac{\partial}{\partial \,x_{\ell}}\,\left(r\circ\Lambda_{i}^{-1}\right) \,\ind_{\Omega_{i}} 
%
\end{equation}
%
%
%
%
\emph {i.e.} the weak derivative does not have pulses/jumps on lower dimensional manifolds.
\begin{proof} Let $\varphi\in C_{0}^{\,\infty}(\Omega)$ then
\begin{multline}\label{Eq integration by parts}
\langle \frac{\partial}{\partial \,x_{\ell}}\,T r, \varphi \rangle_{\scriptscriptstyle D'(\Omega), D(\Omega)}
= - \int_{\Omega} T r\,\frac{\partial \, \varphi}{\partial \,x_{\ell}}
= -  \sum_{i = 1, 2}\int_{\Omega_{i}} r\circ\Lambda^{-1}_{i}\,\frac{\partial \, \varphi}{\partial \,x_{\ell}} \\
%
=  \sum_{i = 1, 2} \int_{\Omega_{i}} \frac{\partial }{\partial \,x_{\ell}}\left(r\circ\Lambda^{-1}_{i}\right)\,\varphi
%
- \sum_{i = 1, 2} \int_{\partial \Omega_{i}} \left(r\circ\Lambda^{-1}_{1}\right)\,\varphi \left(\boldsymbol{\widehat{\nu_{i} } }\cdot\eversor_{\ell}\right) d S . 
%
\end{multline}
We focus on the last two summands. Since $\varphi = 0$ on $\partial\Omega$ this implies
\begin{multline*}
- \sum_{i = 1, 2} \int_{\partial \Omega_{i}} \left(r\circ\Lambda^{-1}_{1}\right)\,\varphi \left(\boldsymbol{\widehat{\nu_{i} } }\cdot\eversor_{\ell}\right) d S\\
= - \int_{\Gamma}\!\! \{\left(r\circ\Lambda^{-1}_{1}\right)\,\varphi \left(\boldsymbol{\widehat{\nu}_{1}  }\cdot\eversor_{\ell}\right)  
%
+ \left(r\circ\Lambda^{-1}_{2}\right)\varphi \left(\boldsymbol{\widehat{\nu}_{2}}\cdot\eversor_{\ell}\right) \}d S = 0 . 
\end{multline*}
The last equality holds since $\boldsymbol{\widehat{\nu}_{1}} = -\boldsymbol{\widehat{\nu}_{2}}$ and $\Lambda_{1}^{-1} = \Lambda_{2}^{-1}$ on $\Gamma$. Combining this fact with \eqref{Eq integration by parts} we conclude \eqref{Eq weak derivatives}.
\end{proof}
\end{lemma}
\begin{theorem}\label{Th change of variable H1 effect}
The map $T$
is an isomorphism from
$H^{1}(\Omega)$ onto itself.
\begin{proof} Since the application $\Lambda$ is a bijection from $\Omega$ into itself the map $T$ is clearly bijective and linear. For the calculation of the norms we use the Change of Variables theorem
\begin{equation*}
\int_{\Omega_{i}} \vert r\circ\Lambda_{i}^{-1}\vert^{2}  =
\int_{\Omega_{i}^{\zeta}} \vert r\vert^{2} \left\vert \det \left(\Lambda_{i}'\right)\right\vert 
\leq 2 \int_{\Omega_{i}^{\zeta}} \vert r\vert^{2}\quad \quad i = 1, 2 ,
\end{equation*}
where the last inequality holds due to \eqref{Eq det Jacobian domain 1}. Equivalently $\Vert
T r\Vert_{0, \Omega_{i}}^{2}\leq 2\,\Vert r\Vert_{0, \Omega_{i}^{\zeta}}^{2}$
i.e. $T$ is a bounded operator in $L^{2}(\Omega)$. For the derivative first consider $u \in C^{\,1}(\Omega)$
take $1\leq \ell\leq N-1$. For the vector function $\Lambda_{\,i}^{-1}:\Omega_{\,i}\rightarrow\Omega_{\,i}^{\,\zeta}$ denote $\Lambda_{ i, k}^{-1}$ its $k$-th component function, thus
\begin{multline*}\label{Ineq first derivatives on isomorphism}
\int_{\Omega_{i}} \left\vert \frac{\partial}{\partial\,x_{\ell}}(u\circ\Lambda_{i}^{-1})\right\vert^{2}  
=
\int_{\Omega_{i}^{\zeta}}  \left\vert \sum_{k\,=\,1}^{N}\frac{\partial u}{\partial\,x_{k}}\,
\frac{\partial \Lambda^{-1}_{i, k}}{\partial\,x_{\ell}}\right\vert^{2} \left\vert
\det \left(\Lambda_{i}'\right)\right\vert \\
= \int_{\Omega_{i}^{\zeta}} \left\vert\frac{\partial u}{\partial\,x_{\ell}}
+ \frac{\partial u}{\partial\,z}\left[1+ (-1)^{i}z\right]\frac{\partial \zeta}{\partial\,x_{\ell}}\right\vert^{2} \left\vert \det \left(\Lambda_{i}'\right)\right\vert\\
\leq 2 \int_{\Omega_{i}^{\zeta}} \left\vert\frac{\partial u}{\partial\,x_{\ell}}\right\vert^{2} 
+ 4 \int_{\Omega_{i}^{\zeta}} \left\vert\frac{\partial u}{\partial\,z}\right\vert^{2}
\left\vert \frac{\partial \zeta}{\partial\,x_{\ell}}\right\vert^{2} 
%
%
\leq \max\,\{2 ,\, 4\,\Vert \zeta \Vert_{\scriptscriptstyle W^{1, \infty}(\Gamma)}^{\,2}\}
\int_{\Omega_{i}^{\zeta}} \left\vert\grad u \right\vert^{2} .
\end{multline*}
For the derivative with respect to $z$ we get
\begin{multline*}
\int_{\Omega_{i}} \left\vert \frac{\partial}{\partial\,z}(u\circ\Lambda_{i}^{-1})\right\vert^{2} 
=
\int_{\Omega_{i}^{\zeta}}  \left\vert \sum_{k\,=\,1}^{N}\frac{\partial u}{\partial\,x_{k}}\,
\frac{\partial \Lambda^{-1}_{i, k}}{\partial z}\right\vert^{2} \left\vert \det \left(\Lambda_{i}'\right)\right\vert \\
= \int_{\Omega_{i}^{\zeta}} \left\vert
\frac{\partial u}{\partial z}\left[1+ (-1)^{i}\zeta\right]\right\vert^{2} \left\vert
\det \left(\Lambda_{i}'\right)\right\vert  
%
%
%
\leq
2\int_{\Omega_{i}^{\zeta}} \left\vert\grad u \right\vert^{2} .
\end{multline*}
Combining both previous inequalities and define
\begin{equation}\label{Def conformal map operator bound}
C_{\,\zeta} \defining  \sqrt{\max\,\{2 ,\, 4\,\Vert \zeta \Vert_{\scriptscriptstyle W^{1, \infty}(\Gamma)}^{\,2}\}} .
\end{equation}
Then it follows 
\begin{equation*}
\Vert \grad (T u) \Vert_{0, \Omega_{i}}^{\,2} 
%
\leq C_{\,\zeta}^{\,2}\,
\Vert \grad u\Vert_{0, \Omega_{i}^{\zeta}}^{\,2}\quad\forall \,u\in C^{\,1}(\Omega) ,
\end{equation*}
for $i = 1, 2$, therefore
\begin{equation*}
\Vert \grad (T u) \Vert_{0, \Omega}
\leq C_{\,\zeta}\,
\Vert \grad u\Vert_{0, \Omega} \quad\forall \,u\in C^{\,1}(\Omega) .
\end{equation*}
The inequality above extends to the whole space $H^{1}(\Omega)$ by density of $C^{\,1}(\Omega)$ in $H^{1}(\Omega)$. Finally, combining the first and second parts we have
\begin{equation}\label{Ineq conformal map operator bound}
\Vert \,T \, r \,\Vert_{1,\Omega}
\leq C_{\zeta}\,
\Vert  \, r \, \Vert_{1,\Omega}\quad\forall \,r\in H^{\,1}(\Omega) .
\end{equation}
i.e. $T$ is a bounded operator on $H^{1}(\Omega)$.
\end{proof}
\end{theorem}
\begin{corollary}\label{Th change of variable Mv effect}
The map $T$ is an isomorphism from $V$ onto itself.
\begin{proof} Observe that $\Lambda_{1}^{-1}\vert_{\partial\Omega_{1} - \Gamma} = I \, \vert_{\partial\Omega_{1} - \Gamma}$, then $T r = 0$ on $\partial\Omega_{1} - \Gamma$ i.e. $T$ is a bijection from $V$ into itself and due to previous theorem the result follows.
\end{proof}
\end{corollary}
%
%
%
%
%
%
%
%
\subsection{The Fractional Mapping Operator on $\pert(\Omega, \Gamma)$}  
Consider the application $T:\apert(\Gamma, \Omega)\rightarrow \mathcal{L}(H^{1}(\Omega))$, where $\zeta\mapsto T(\zeta)$ is defined by equation \eqref{Eq morphism of functions}. Since $\apert(\Gamma, \Omega)$ is not a linear space, only a convex set, $T$ can not be linear; however $T$ does not respect convex combinations either. Therefore, the nonlinearity of $T$ does not lie only on its domain of definition but also on its algebraic structure. Clearly $T(0) = I$, we will show $T$ is continuous at $0$ in the pointwise topology.
\begin{lemma}\label{Th continuity of the conformal map operator on smooth functions}
Let $\{\zeta_{\,n}\}\subseteq \apert(\Gamma, \Omega)$ bounded in $ W^{1, \infty}(\Gamma)$ and such that $\Vert \zeta_{\,n}\Vert_{\scriptscriptstyle  C(\Gamma)}\rightarrow 0$, then
\begin{equation*}
\left\Vert T(\zeta_{\,n}) u -  u\right\Vert _{H^{1}(\Omega)}\rightarrow 0\,,\quad\forall\, u\in C^{1}(\Omega) .
\end{equation*}
%
%
%
%
\begin{proof}
First notice that for $i = 1, 2$ 
\begin{equation*}
\Lambda_{i}^{-1}(\zeta_{\,n})(\Xthilde, X_{\scriptscriptstyle N}) = 
(\xthilde, z (1 - (-1)^{i})\zeta_{n}(\xthilde) + \zeta_{n}(\xthilde)) .
\end{equation*}
Then the uniform convergence $\Vert \zeta_{\,n} \Vert_{\scriptscriptstyle C(\Gamma) }\rightarrow 0$ implies 
\begin{equation}\label{Eq conformal map convergence}
\Vert \Lambda_{1}^{-1}(\zeta_{\,n})\ind_{\scriptscriptstyle \Omega_{1}}
+ I\ind_{\scriptscriptstyle \Gamma}
+ \Lambda_{2}^{-1}(\zeta_{\,n}) \ind_{\scriptscriptstyle \Omega_{2}}
- I\Vert_{\scriptscriptstyle C(\Omega)}\xrightarrow[n\rightarrow \infty]{} 0 .
\end{equation}
Here $I\,\ind_{\scriptscriptstyle \Gamma}$ has to be introduced for the convergence in the space of continuous functions; also recall the fact that $\Lambda_{1}^{-1}(\zeta) \vert_{\Gamma}= \Lambda_{2}^{-1}(\zeta)\vert_{\Gamma} = I\vert_{\Gamma}$ for all $\zeta\in \apert(\Gamma, \Omega)$. Take $u\in C^{1}(\Omega)$ and $\x\in \Omega$ fixed, then due to the Mean Value Theorem in multiple dimensions we have
\begin{equation*}
\vert T(\zeta_{\,n}) u (\x) - u(\x) \vert
= \vert u\circ \Lambda_{i}^{-1}(\zeta_{\,n}) (\x) - u (\x) \vert \\
\leq \vert \grad u (\boldsymbol{\xi} )  \vert\, \vert\Lambda_{i}^{-1}(\zeta_{\,n})(\x) - \x \vert 
%
\end{equation*}
Where $\boldsymbol {\xi} $ lies in the segment uniting $\x$ and $\Lambda_{i}^{-1}(\zeta_{\,n})(\x)$. Thus 
%
%
\begin{equation}\label{Ineq norm convergence map operator estimate}
\int_{\Omega_{i}} \vert T(\zeta_{\,n}) u - u \vert^{\,2} 
\leq  \Vert \,u \,\Vert_{\scriptscriptstyle H^{1}(\Omega)}^{2} \,
\Vert \Lambda_{i}^{-1}(\zeta_{\,n}) - I\Vert_{\scriptscriptstyle C(\Omega_{i})}^{2}
\quad \text{for}\; i = 1, 2.
\end{equation}
Due to \eqref{Eq conformal map convergence} we conclude
\begin{equation}\label{Eq L2 convergence conformal map operator}
\Vert T(\zeta_{\,n}) u - u\Vert_{\scriptscriptstyle 0,\Omega}\rightarrow 0 .
\end{equation}
For the $H^{1}(\Omega)$-convergence first notice that due to the inequality \eqref{Ineq conformal map operator bound} and definition \eqref{Def conformal map operator bound} it follows 
\begin{equation}
\Vert T(\zeta_{\,n}) u\Vert_{\scriptscriptstyle 1,\Omega}\leq
\sup _{n\,\in\,\N}\sqrt{\max \{2 , 4 \Vert \zeta_{\,n} \Vert_{\scriptscriptstyle W^{1, \infty}(\Omega)}^{2}\}} \; 
\Vert u\Vert_{\scriptscriptstyle 1, \Omega} .
\end{equation}
Where the supremum is finite because of the boundedness of $\{\zeta_{\,n}\}$ in $W^{1, \infty}(\Gamma)$. Then the sequence $\{T(\zeta_{\,n})u\}$ has a weakly convergent subsequence in $H^{1}(\Omega)$ and due to \eqref{Eq L2 convergence conformal map operator} the weak limit must be $u$. Moreover, the Rellich-Kondrasov compactness theorem 
implies that the whole sequence converges weakly to the same limit $u$. 

Next we prove that the $H^{1}(\Omega)$-norms converge. The strong convergence in $L^{2}(\Omega)$ is given by the Rellich-Kondrachov theorem, therefore we focus only on the derivatives. For any $1\leq \ell\leq N-1$ we have
\begin{multline*}
\int_{\Omega_{i}} \left\vert \frac{\partial}{\partial\,x_{\ell}}\left[u\circ\Lambda_{i}^{-1}(\zeta_{\,n})\right]\right\vert^{2}  
= \int_{\Omega_{i}^{\zeta_{n}}} \vert\partial_{\ell} u
+ \partial_{z} u(1+ (-1)^{i}z)
\partial_{\ell} \zeta_{\,n}\vert^{2} \vert
\det \Lambda_{i}'(\zeta_{\,n})\vert\\
= \int_{\Omega_{i}^{\zeta_{n}}} \vert\partial_{\ell} u 
+ \partial_{z} u(1+ (-1)^{i}z)
\partial_{\ell} \zeta_{\,n} \vert^{2} \vert
\det \Lambda_{i}'(\zeta_{\,n})\vert\\
= \int_{\Omega_{i}^{\zeta_{n}}} \vert\partial_{\ell} u
+ \partial_{z} u (1+ (-1)^{i}z)
\partial_{\ell} \zeta_{\,n}\vert^{2}
(1-(-1)^{i}\zeta_{\,n} ) .
\end{multline*}
On one hand, it is clear that the integrand converges to $\vert \dfrac{\partial \, u}{\partial\,x_{\ell}} \vert^{2}\,\ind_{\Omega_{i}}$, on the other hand, we have the estimate
\begin{multline*}
\left\vert \partial_{\ell}  u
+ \partial_{z}  u(1+ (-1)^{i}z) \partial_{\ell} \zeta_{\,n}\vert^{2}
(1-(-1)^{i}\zeta_{\,n}) \ind_{\Omega_{i}^{\zeta_{n} } }\right\vert \\
\leq 2 \, \{\vert\partial_{\ell} u 
\vert^{2} +
\vert 1+ (-1)^{i}z \vert
\vert \partial_{\ell} \zeta_{\,n}  \vert^{2}
 \vert \partial_{z} u \vert^{2} \} \, 2 
\leq 4 \, \{\vert\partial_{\ell} \, u
\vert^{2} +
2 \,\Vert\zeta_{\,n}\Vert^{2}_{\scriptscriptstyle W^{1, \infty}(\Gamma)}
\vert \partial_{z} \, u  \vert^{2} \}\\
\leq 4 \max\{1 + 2\,\sup_{n\,\in\,\N}\Vert\zeta_{\,n}\Vert^{2}_{\scriptscriptstyle W^{1, \infty}(\Gamma)}\}\,
\vert\grad u\,\vert^{ 2} \in L^{1}(\Omega)
\end{multline*}
Thus, Lebesgue's dominated convergence theorem yields
%
%
\begin{equation}\label{Eq convergence of l-th derivative}
\left\Vert \frac{\partial }{\partial\,x_{\ell}}\,(T(\zeta_{\,n})u)\right\Vert^{2}_{0, \Omega_{i}}
\xrightarrow [n\,\rightarrow\,\infty]{}
\left\Vert \frac{\partial u}{\partial\,x_{\ell}} \right\Vert^{2}_{0, \Omega_{i}} 
 1\leq \ell\leq N - 1,\;\text{for}\; i = 1, 2 .
\end{equation}
%
%
For the derivative with respect to $z$ we get
\begin{multline*}
\int_{\Omega_{i}} \left\vert \frac{\partial}{\partial\,z}\left[u\circ\Lambda_{i}^{-1}(\zeta_{\,n})\right]\right\vert^{2} 
= \int_{\Omega_{i}^{\zeta_{n}}} \left\vert
\frac{\partial u}{\partial z}\left[1+ (-1)^{i}\zeta_{\,n}\right]\right\vert^{2} \left\vert \det \left(\Lambda_{i}'(\zeta_{\,n})\right)\right\vert\\
= \int_{\Omega} \left\vert
\frac{\partial u}{\partial z}\right\vert^{2}
\left[1+ (-1)^{i}\zeta_{\,n}\right]^{3}\,\ind_{\Omega_{i}^{\zeta_{n}}}
\end{multline*}
Again, the integrand converges to $\vert \dfrac{\partial u}{\partial z}\vert^{2}\,\ind_{\Omega_{i}}$ and it is bounded by $2 \,\vert \dfrac{\partial u}{\partial z}\vert^{\,2}$ which is an element of $ L^{1}(\Omega)$. Hence, the Lebesgue dominated convergence theorem yields
\begin{equation}\label{Eq convergence of z derivative}
\left\Vert \frac{\partial }{\partial z}\,(T(\zeta_{\,n})u)\right\Vert^{2}_{0,\,\Omega_{i}}
\xrightarrow [n\,\rightarrow\,\infty]{}
\left\Vert \frac{\partial \, u}{\partial z} \right\Vert^{2}_{0,\,\Omega_{i}}
, \; \text{for}\;  i = 1, 2 .
\end{equation}
The equations \eqref{Eq convergence of l-th derivative} and \eqref{Eq convergence of z derivative} give the convergence of the $L^{2}(\Omega)$-norms of the gradients $\Vert \grad T(\zeta_{\,n})u\Vert_{\,0,\Omega} \rightarrow \Vert \grad u\Vert_{\,0,\Omega}$ and then $\Vert  T(\zeta_{\,n})u\Vert_{\,1,\Omega} \rightarrow \Vert  u\Vert_{\,1,\Omega}$. This fact together with the weak convergence in $H^{1}(\Omega)$ imply $\Vert  T(\zeta_{\,n})u - u \Vert_{\,1,\Omega} \rightarrow 0$.
\end{proof}
\end{lemma}
\begin{theorem}\label{Th continuity of the conformal map operator}
Let $\{\zeta_{\,n}\}\subseteq \apert(\Gamma, \Omega)$ bounded in $ W^{1, \infty}(\Gamma)$ and such that $\Vert \zeta_{\,n}\Vert_{\scriptscriptstyle  C(\Gamma)}\rightarrow 0$, then
\begin{equation*}
\left\Vert T(\zeta_{\,n})\, r -  r\right\Vert _{\scriptscriptstyle H^{1}(\Omega)}\rightarrow 0\,,\quad\forall\, r\in H^{1}(\Omega)
\end{equation*}
\emph{i.e.} $T(\zeta_{\,n})$ converges to the identity $I$ in the strong operator topology.
\begin{proof}
We use the standard density argument. Let $r\in H^{1}(\Omega)$, take $\{ u_{\scriptscriptstyle j}\}\subseteq C^{\,1}(\Omega)$ such that $\Vert  u_{\scriptscriptstyle j} - r\Vert_{1,\Omega}\rightarrow 0$. Recall definition \eqref{Def conformal map operator bound} and inequality \eqref{Ineq conformal map operator bound}, then
\begin{multline*}
\Vert \,T(\zeta_{\,n}) \, r -  r\Vert _{\scriptscriptstyle H^{1}(\Omega)} 
\leq \Vert \,T(\zeta_{\,n}) \, r -  T(\zeta_{\,n})  u_{\scriptscriptstyle j} \Vert _{\scriptscriptstyle H^{1}(\Omega)} 
+ \Vert \,T(\zeta_{\,n})  u_{\scriptscriptstyle j} -   u_{\scriptscriptstyle j}\Vert _{\scriptscriptstyle H^{1}(\Omega)}
+ \Vert \, u_{\scriptscriptstyle j} -  r\Vert _{\scriptscriptstyle H^{1}(\Omega)}\\[3pt]
\leq ( 1 + \Vert \,T(\zeta_{\,n})\Vert)  \Vert  u_{\scriptscriptstyle j} -  r\Vert _{\scriptscriptstyle H^{1}(\Omega)} 
+ \Vert \,T(\zeta_{\,n})  u_{\scriptscriptstyle j} -   u_{\scriptscriptstyle j} \Vert _{\scriptscriptstyle H^{1}(\Omega)}\\
\leq  ( 1 + 
\sup_{k\,\in\,\N} \sqrt{\max\,\{2, \,4\,\Vert\zeta_{\,k}\Vert_{\scriptscriptstyle W^{1, \infty}(\Gamma)}\}}) \; 
 \Vert \, u_{\scriptscriptstyle j} -  r \Vert _{\scriptscriptstyle H^{1}(\Omega)} 
+
 \Vert \,T(\zeta_{\,n})  u_{\scriptscriptstyle j} -   u_{\scriptscriptstyle j} \Vert _{\scriptscriptstyle H^{1}(\Omega)}
\end{multline*}
Fix $j\in \N$ such that the first summand of the right hand side is less than $\tfrac{\epsilon}{2}$. Due to lemma \ref{Th continuity of the conformal map operator on smooth functions} there exists $N\in \N$ such that $n\geq N$ implies the second summand on the right hand side of the expression above is less than $\tfrac{\epsilon}{2}$ and the proof is complete.
\end{proof}
\end{theorem}
\begin{corollary}\label{Th continuity of the conformal map operator inverse}
Let $\{\zeta_{\,n}\}\subseteq \apert(\Gamma, \Omega)$  bounded in $ W^{1, \infty}(\Gamma)$ and such that $\Vert \zeta_{\,n}\Vert_{\scriptscriptstyle  C(\Gamma)}\rightarrow 0$, then
\begin{equation*}
\left\Vert \, T^{\,-1}(\zeta_{\,n})\,r - r \, \right\Vert _{ \scriptscriptstyle H^{1}(\Omega)}\xrightarrow[n\rightarrow \infty]{} 0
\end{equation*}
\emph{i.e.} the sequence of inverse operators converge.
\begin{proof} Repeating the techniques exposed in Lemma \ref{Th continuity of the conformal map operator on smooth functions} and in theorem \ref{Th continuity of the conformal map operator} applied to the fractional bijective maps $\{\Lambda_{i}(\zeta_{\,n})\}$ defined in \eqref{Def conformal map domain i} the result follows.
\end{proof}
\end{corollary}
%
%
%
%
\begin{remark}
Using theorem \ref {Th continuity of the conformal map operator} it can be proved that $\left\Vert \, T(\zeta_{\,n}) - \jmath \, 
\right\Vert _{\scriptscriptstyle \mathcal{L}(C^{\,1}(\Omega), \,L^{2}(\Omega))}\rightarrow 0$ as $\Vert \zeta_{n}\Vert_{\scriptscriptstyle C(\Gamma)} \rightarrow 0$ if $\{\zeta_{n}\}$ is bounded in $W^{1, \infty}(\Gamma)$. Here $\jmath:C^{1}(\Omega) \hookrightarrow L^{2}(\Omega)$ is the embedding operator $\jmath (\varphi) \defining \varphi$ and $T(\zeta_{\,n})$ is regarded as an operator from $C^{\,1}(\Omega)$ to $L^{2}(\Omega)$.
Also, using the same technique in obtaining the estimate \eqref{Ineq norm convergence map operator estimate} we can show $\left\Vert \, T(\zeta_{\,n}) - \jmath \, \right\Vert _{\scriptscriptstyle \mathcal{L}(C^{2}(\Omega), \,H^{1}(\Omega))}\rightarrow 0$ with $\jmath:C^{2}(\Omega) \hookrightarrow H^{1}(\Omega)$; i.e. in order to get convergence in the norm of the operators higher degrees of regularity are needed. 
\end{remark}
%
%
%
%
\subsection{The Flattened Problem}  
Consider the problem \eqref{Pblm Direct Formulation Perturbed} subject to the changes of variable described above. Introducing \eqref{Eq det Jacobian domain 1} and \eqref{Eq gradient structure domain 1} in each summand of left hand side in \eqref{Pblm Direct Formulation Perturbed} we have
\begin{equation*}
\int_{\Omega_{i}^{\zeta}}\grad_{\!\! \scriptscriptstyle X}\,q\cdot\grad_{\!\! \scriptscriptstyle X} \,r =
\int_{\Omega_{i}} \grad_{\!\! \scriptscriptstyle X}\,T q\cdot\grad_{\!\! \scriptscriptstyle X}  T r
\,\left\vert \det \Lambda_{i}'\right\vert \\
=\int_{\Omega_{i} } (A_{i}^{\zeta})^{T} A_{i}^{\zeta}\,\grad_{\!\! x} T q 
\cdot\grad_{\!\! x} T r \,(1 - (-1)^{i}\zeta). 
%
\end{equation*}
%
%
%
%
With
\begin{equation}\label{Eq matrix inner product}
 (A_{i}^{\zeta})^{T} A_{i}^{\zeta}\\
 = \left[ \!\!\! \begin{array}{cc}
 (1- (-1)^{i}\zeta)I & \!\!\! (-1)^{i} (z - (-1)^{i} )\gradth \zeta\\[10pt]
 (-1)^{i}(z - (-1)^{i} ) \gradth^{T} \zeta &  \!\!\! \dfrac{\vert(z - (-1)^{i} )\grad\zeta\vert^{2}+1}{1 - (-1)^{i}\zeta}
 \!\!\!
\end{array}
\right] .
\end{equation}
Due to corollary \ref{Th change of variable Mv effect} the quantifiers $\forall\;T r\in V$ and $\forall\,r\in V$ are equivalent, therefore we conclude that the solution $q$ of problem \eqref{Pblm Direct Formulation Perturbed} satisfies the variational statement
%
%
\begin{multline}\label{Eq flattened perturbed magnified problem}
q\in V:
\sum_{i \, = \, 1, \,2}\int_{\Omega_{i}}
\frac{1-(-1)^{i}\,\zeta }{\epsilon^{\,i-1}} \, (A_{i}^{\zeta})^{T} A_{i}^{\zeta} \,\grad T q \cdot \grad r\\
%
%
=\int_{\Gamma} \frac{1}{\vert(-\gradth\zeta, 1)\vert}\left(T f\right) \, r\;d\,S 
+ \sum_{i \, = \, 1, \,2} \,\int_{\Omega_{i}}(1-(-1)^{i}\,\zeta\,)\left(T F\right) r
%
\, , \quad\forall\;r\in V .
\end{multline}
%
%
In the first summand of the left hand side, the notation $Tf$ stands for $f\circ\Lambda_{1}^{-1}$ or $f\circ\Lambda_{2}^{-1}$ indistinctively since $\Lambda_{1}^{-1} = \Lambda_{2}^{-1}$ on $\Gamma$. The surface integral summand implicitly uses the fact that the upwards unitary vector normal to the surface $\Gammaz$ is given by $\n = \dfrac{(-\gradth\zeta, 1)}{\vert(-\gradth\zeta, 1)\vert}$. Declaring $\varrho \defining T q$ the problem \eqref{Eq flattened perturbed magnified problem} is equivalent to
%
%
\begin{multline}\label{Eq geometry effect on inner product problem}
\varrho\in V:
\sum_{i \, = \, 1, \,2 }\int_{\Omega_{i}}
\frac{ 1-(-1)^{i}\,\zeta}{\epsilon^{\,i-1}} \, (A_{i}^{\zeta})^{T} A_{i}^{\zeta} \;\grad \varrho \cdot \grad r\\
=\int_{\Gamma} \frac{1}{\vert(-\gradth\zeta, 1)\vert}\left(T f\right) \, r\;d\,S 
+ \sum_{i \, = \, 1, \,2 } \int_{\Omega_{i}}(1-(-1)^{i}\,\zeta\,)\left(T F\right) r
\, , \quad\forall\;r\in V .
\end{multline}
%
%
Next, we focus on some properties of the involved matrices.
\begin{lemma}\label{Th inverse matrix characterization}
For $i = 1, 2$ the inverse matrix of $A_{i}^{\zeta}$ is given by
\begin{equation}\label{Eq inverse matrix characterization}
(A_{i}^{\zeta})^{-1} = \left[\begin{array}{cc}
I & (1 - (-1)^{i} z)\gradth \zeta \\[5pt]
\0^{T} & 1 - (-1)^{i} \zeta
\end{array}
\right]
\end{equation}
\begin{proof}
By direct calculation 
\end{proof}
\end{lemma}
\begin{corollary}\label{Th norm of the conformal map}
Let $(\xthilde, z)\in \Omega_{i}^{\zeta}$ be fixed with $i \in\{  1, 2\}$. Then the linear operator $\boldsymbol{\xi}\mapsto (A_{i}^{\zeta}(\x))^{-1} \boldsymbol{\xi}$ from $\R^{\! N}$ into itself endowed with the canonical inner product satisfies
\begin{equation}\label{Eq norm of the conformal map}
\Vert (A_{i}^{\zeta}(\x) )^{-1}\Vert_{\scriptscriptstyle\mathcal{L}(\R^{\!N} ) } \leq 
\sqrt{N + 3 + 4\,\Vert \zeta\Vert_{\scriptscriptstyle W^{1, \infty}(\Gamma)}^{2}} \\
\, , \;\; \text{for all}\;\,\x\in \Omega_{i}^{\zeta},\;\,\text{and}\;\,i \in\{1,2\}
\end{equation}
\begin{proof}
We compute the Frobenius norm of the operator adding the squared inner product norms $\vert \cdot \vert$ of each column vector in \eqref{Eq inverse matrix characterization}. This gives
\begin{equation*}
\Vert (A_{i}^{\zeta}(\x))^{-1}\Vert_{\scriptscriptstyle \mathcal{L}(\R^{\! N} ) }^{2}
\leq (N-1) \\
+ (1 - (-1)^{i}z)^{2} 
\vert \grad \zeta (\xthilde)\,\vert^{2} + (1 - (-1)^{i} \zeta )^{2} .
\end{equation*}
Since $z, \zeta(\xthilde)\in [-1,1]$ for all $\x = (\xthilde, z)\in \Omega$ we estimate $\vert 1\pm z\vert\leq 2$ and $\vert 1 \pm\zeta\vert \leq 2$. The gradient of $\zeta$ is estimated by the $W^{1, \infty}(\Gamma)$-norm and the result follows.
\end{proof}
\end{corollary}
\begin{proposition}\label{Th uniform coercivity of matrices}
The matrices $(1 - (-1)^{i}\zeta) (A_{i}^{\zeta})^{T} A_{i}^{\zeta}$ are uniformly coercive in $\R^{\! N}$ for $i = 1, 2$ \emph{i.e.} there exists $e(\zeta) > 0$ such that
\begin{equation}\label{Eq uniform coercivity of matrices}
e(\zeta)\,\vert \xiv\vert^{2} \leq (1 - (-1)^{i}\zeta) (A_{i}^{\zeta})^{T} A_{i}^{\zeta} \,\xiv\cdot\xiv
%
\end{equation}
For all $\xiv\in \R^{\! N}$ and each $(\xthilde, z)\in \Omega$. 
\begin{proof}
Let $\x\in \Omega_{i}$ be fixed and $\xiv\in \R^{N}$ unitary then
\begin{multline*}
(1 - (-1)^{i}\zeta) (A_{i}^{\zeta})^{T} A_{i}^{\zeta} \,\xiv\cdot\xiv 
= (1 - (-1)^{i}\zeta)  A_{i}^{\zeta} \,\xiv\cdot A_{i}^{\zeta}\, \xiv
\geq (1 - \Vert \zeta\Vert_{\scriptscriptstyle C(\Gamma)}) 
\left\vert A_{i}^{\zeta} \,\xiv \right\vert^{2} \\
\geq (1 - \Vert \zeta\Vert_{\scriptscriptstyle C(\Gamma)})\, \min_{\vert \etav\vert = 1} \left \vert A_{i}^{\zeta} \,\etav\right\vert^{2}
= \frac{1 - \Vert \zeta\Vert_{\scriptscriptstyle C(\Gamma)}}{\Vert (A_{i}^{\zeta})^{-1}\Vert_{\scriptscriptstyle \mathcal{L}(\R^{\! N})}^{2}} 
\geq \frac{1 - \Vert \zeta\Vert _{C(\Gamma)}}{N + 3 + 4\,\Vert \zeta \Vert^{2}_{\scriptscriptstyle W^{1, \infty}(\Gamma)}} .
\end{multline*}
Where, the last inequality comes from corollary \ref{Th norm of the conformal map}. Defining $e(\zeta)\defining \frac{1 - \Vert \zeta\Vert _{\scriptscriptstyle C(\Gamma)}}{N + 3 + 4\,\Vert \zeta \Vert^{2}_{\scriptscriptstyle W^{1,\infty}(\Gamma)}}$ the statement \eqref{Eq uniform coercivity of matrices} holds.
\end{proof}
\end{proposition}
\begin{corollary}\label{Th equivalent inner product}
The form $[\,\cdot, \,\cdot]:V\times V\rightarrow \R$
\begin{equation}\label{Def bilinear form}
[\kappa, \pi] \defining \sum_{i \,=\,1,\,2} \int_{\Omega_{i}} \!\! \frac{1 - (-1)^{i}\zeta}{\epsilon^{\,i - 1}} \, 
(A_{i}^{\zeta})^{T} A_{i}^{\zeta}\, \grad \kappa\cdot\grad \pi 
\end{equation}
defines an inner product on $V$ which induces the same topology as the induced by the standard inner product on $H^{1}(\Omega)$.
\begin{proof}
The form \eqref{Def bilinear form} is well-defined since the application
\begin{equation*}
(\xthilde, z)\in\Omega_{i}\mapsto (1 - (-1)^{i}\zeta(\xthilde))
(A_{i}^{\zeta}(\xthilde, z))^{T} A_{i}^{\zeta}(\xthilde, z)
\end{equation*}
for $i = 1, 2$ is in $L^{\infty}(\Omega_{\,i}, \R^{\! N\times N})$. Clearly the form is bilinear and symmetric. For the continuity we have $[\kappa, \pi] \leq \tfrac{2}{\epsilon}\,\Vert A_{i}^{\zeta}\Vert_{\infty}^{2}\,\Vert \grad \kappa\Vert_{ 0,\,\Omega}\,\Vert \grad \pi\Vert_{ 0,\,\Omega}$ $\leq \tfrac{2}{\epsilon}\,\Vert A_{i}^{\zeta}\Vert_{\infty}^{2}\,\Vert \kappa\,\Vert_{ 1,\,\Omega}\,\Vert \pi\Vert_{ 1,\,\Omega}$ for all $\kappa, \pi\in V$; in particular $[\kappa, \kappa]\leq \tfrac{2}{\epsilon}\,\Vert A_{i}^{\zeta}\Vert_{\infty}^{2}\,\Vert \kappa\,\Vert_{ 1,\,\Omega}^{2}$. The homogeneous condition of the induced norm comes from the uniform coercivity of the matrices shown in proposition \ref{Th uniform coercivity of matrices}. Hence
\begin{equation*}
\frac{e(\zeta)}{1 + C_{\scriptscriptstyle \Omega}}\,\Vert\kappa\Vert_{ 1,\Omega}^{2} \leq e(\zeta)\,\Vert\grad \kappa\Vert_{ 0,\,\Omega}^{\,2}\leq [\kappa, \kappa]\quad\forall\,\kappa\in V .
\end{equation*}
Where $C_{\scriptscriptstyle\Omega}$ is the Poincar\'e constant associated to the domain $\Omega$ valid for all elements of $V$ given the boundary conditions. Therefore, the induced norms are equivalent.
\end{proof}
\end{corollary}
\begin{theorem}\label{Th conformally mapped well-posedness}
The problem \eqref{Eq geometry effect on inner product problem} is well-posed.
%
%
\begin{proof}
The equivalence of norm induced by the inner product $[\,\cdot, \cdot]$ to the standard $H^{1}(\Omega)$-norm on $V$ is shown in corollary \ref{Th equivalent inner product}, therefore the well-posedness of problem \eqref{Eq geometry effect on inner product problem} follows from Lax-Milgram's lemma.
\end{proof}
\end{theorem}
%
%
%
%
%
%
\subsection{Geometric Perturbation and Inner Products}\label{Sec Geometric Perturbation and Inner Products}
This section is aimed to analyze the highly nonlinear impact of the geometry in terms of the inner product. Consider the bounds
\begin{equation}\label{Stmt global convergence}
\Vert p - q^{\zeta} \Vert_{\scriptscriptstyle V} \leq 
\Vert p - T^{-1}(\zeta) p \, \Vert_{\scriptscriptstyle V} + 
\Vert T^{-1}(\zeta) p  - q^{\zeta} \Vert_{\scriptscriptstyle V}
\end{equation}
The first summand converges due to theorem \ref{Th continuity of the conformal map operator}. For the convergence of the second summand, due to corollary \ref{Th continuity of the conformal map operator inverse} it is equivalent to prove
\begin{equation}\label{Stmt convergence inner products}
\Vert p  - T(\zeta) q^{\zeta} \Vert_{\scriptscriptstyle V} = 
\Vert p  - \varrho^{\zeta} \Vert_{\scriptscriptstyle V}\rightarrow 0\;
\text{as} \; \Vert \zeta\Vert_{\scriptscriptstyle C(\Gamma)}\rightarrow 0 \;
\text{with} \; \Vert \zeta\Vert_{\scriptscriptstyle W^{1, \infty}(\Gamma)}\; \text{bounded} .
\end{equation}
We are to estimate the norm above by comparing the problems \eqref{Pblm Direct Formulation Original} and \eqref{Eq geometry effect on inner product problem}. In problem \eqref{Eq geometry effect on inner product problem} the effect of the geometry on the inner product structure is fully contained in the matrices 
\begin{equation}\label{Def family of matrices}
\{(1 - (-1)^{i} \zeta) (A^{\zeta}_{i})^{T} A^{\zeta}_{i}:\x\in \Omega_{i}\}\quad i = 1, 2.
\end{equation}
Notice this is a family of symmetric and therefore diagonalizable matrices. However, they depend on the point $\x\in \Omega$ and, in general, they do not commute for $\x, \x '\in \Omega$ different. Therefore, it can not be assured that the family \eqref{Def family of matrices} is simultaneously diagonalizable. Observe that the matrices \eqref{Eq extended gradient structure domain 1} have entries multiplied by the factors $\dfrac{1}{1 - (-1)^{i} \zeta(\xthilde)}$, for $ i = 1, 2$ respectively. This implies that the maps induced by the matrices are not linear with respect to the perturbation $\zeta$. The following hypothetic assumptions illustrate the nonlinearity of the dependence.
\begin{enumerate}[(i)]
\item
Suppose that $\zeta$ is a piecewise linear affine function. Although $\gradth\zeta$ is piecewise constant, the map $\x\mapsto (1 - (-1)^{i} \zeta) (A^{\zeta}_{i})^{T} A^{\zeta}_{i}$ is not piecewise constant.  

\item
Assume that the family of matrices \eqref{Def family of matrices} is diagonal. Testing the problems \eqref{Pblm Direct Formulation Original} and \eqref{Eq geometry effect on inner product problem} 
with $p - \varrho$ and subtracting them yields
%
%
\begin{multline}\label{Eq comparison with eigenvalues}
\sum_{i\,=\,1, \,2}\int_{\Omega_{i}} \frac{\grad p - \grad \,\varrho}{\epsilon^{i-1}} \cdot
\{\,
\grad p - \left[\begin{array}{ccc}
\mu_{1} &  \ldots & 0\\
\vdots  & \vdots & \vdots \\
0 & \ldots & \mu{\scriptscriptstyle N}
\end{array}\right]\grad \,\varrho\}\\
= \int_{\Omega} F(p -  \varrho) - \sum_{i\,=\,1, \,2}\int_{\Omega_{i}} (1- (-1)^{i}\,\zeta) (T F)(p -  \varrho)\,d S\\
+\int_{\Gamma}f (p -  \varrho)\,d S
- \int_{\Gamma} \frac{1}{\vert (-\gradth \zeta, \,1)\vert}\,(T f)\,(p -  \varrho)\,d S
\end{multline}
%
%
Where $\mu_{1}, \ldots, \mu_{\scriptscriptstyle N}$ are the eigenvalues. Nevertheless $\mu_{j} = \mu_{j}(\x)$ for $1\leq j\leq N$ i.e. they depend on the position $\x$ within the domain $\Omega$. 
\end{enumerate}
We close this section reviewing the simplest possible scenario.
%
%
%
%
%
%
\subsection*{Fractional Mapping of Domains in the One Dimensional Case}
In this section we analyze the fractional mapping technique in the one dimensional setting. For this case the problem \eqref{Eq geometry effect on inner product problem} reduces to
%
%
\begin{multline}\label{Pblm conformally mapped problem 1-d}
\frac{1}{1 + \zeta}\int_{-1}^{\,0}\del \varrho \cdot \del r
+ \frac{1}{\epsilon}\,\frac{1}{1 - \zeta}\int_{0}^{\,1}\del \varrho \cdot \del r
= f(\zeta)\, r (0) \\
+\sum_{i = 1, 2} (1 - (-1)^{i} \zeta)\int_{\frac{(-1)^{i} - 1}{2} }^{\frac{(-1)^{i} + 1}{2}} r(x) \,F\left(\frac{x - \zeta}{1 - (-1)^{i} \zeta}\right) \,dx
%
\end{multline}
%
%
In order to get a-priori estimates the left hand side of the problem \eqref{Pblm conformally mapped problem 1-d} suggests testing equations \eqref{Eq direct original 1-d} and \eqref{Pblm conformally mapped problem 1-d} with the function
%
%
\begin{equation}\label{Def piecewise linear combination}
( p - \frac{1}{1 + \zeta}\,\varrho)\,\ind_{(-1, \,0]} +
( p - \frac{1}{1 - \zeta}\,\varrho\,)\ind_{[0, \,1)} 
= p
- \{\frac{1}{1 + \zeta}\,\varrho\,\ind_{(-1, \,0]} + \frac{1 }{1 - \zeta}\,\varrho\,\ind_{[0,\,1)}\} .
\end{equation}
%
%
However, the test function presented above (as the second summand in the right hand side illustrates) is not eligible, because it does not belong to $V$ unless $\zeta = 0$ or $\varrho (0) = 0$. The first condition removes the perturbation leaving the original problem \eqref{Eq direct original 1-d} and the second can not be assured. Any other attempt of estimating $\Vert p - \varrho\Vert_{\scriptscriptstyle V}$ demands test functions equivalent to \eqref{Def piecewise linear combination}, because of the coefficients disagreement in problems \eqref{Eq direct original 1-d} and \eqref{Pblm conformally mapped problem 1-d}. Of course such piecewise functions are not in the test space $V$ due to the trace continuity requirements.

In the one dimensional case, the fractional mapping technique defines an inner product much easier to understand than the corresponding to the multidimensional case \eqref{Eq geometry effect on inner product problem}. It is clear that the dependence of the inner product with respect to the perturbation is piecewise linear fractional as the map $\Lambda: (0,1)\rightarrow (0,1)$ itself. Additionally, direct calculations can be done to find explicitly, the dependence of the eigenvalues associated to the problem \eqref{Eq direct perturbed 1-d}. This dependence also turns out to be piecewise fractional on $\zeta$, closely related to $\Lambda$. For a deep discussion on boundary perturbation of the Laplace eigenvalues see \cite{StrangGrinfeld2004,  StrangGrinfeld2012}. 

Finally,  the question of strong convergence \eqref {Stmt convergence inner products} can be solved using the Rellich-Kondrachov theorem. However, this is not a constructive result and it does not yield explicit estimates such as inequality \eqref{Ineq continuous dependence wrt left perturbation} in section \ref{Sec global estimate of perturbation}. 
%
%
%
%
%
%
%
%
%
%
%
%
%
%
\section{Concluding Remarks and Discussion}
The present work yields several conclusions and open questions. 
\begin{asparaenum}[(i)]
\item
In section \ref{Sec stability estimates}, extra hypothesis of regularity on the forcing terms involved were introduced in order to conclude weak convergence in a first step, and strong convergence in a second one. Also, for the geometric perturbations $\zeta \in \pert(\Omega, \Gamma)$ it is not enough to have convergence in $C(\Gamma)$, there is also need for boundedness in $W^{1, \infty}(\Gamma)$. 

\item
The conditions of convergence for the interface are acceptable in the context of saturated porous media fluid flow.  Moreover, for the modeling of saturated fluid flow through deformable porous media in the \emph{elastic regime}, these conditions are natural, because for high gradients of deformation the \emph{elasto-plastic} and \emph{plastic} regimes start taking place, see \cite{Morales1}. 
\item
Mapping the perturbed domains with fractional applications as in section \ref{Sec comparing problems}, decomposes the nonlinearity of the question in two parts: the fractional mapping operator $T: \pert(\Omega, \Gamma)\rightarrow \mathcal{L}(H^{1}(\Omega))$ which is clearly nonlinear, and the effect of the geometric perturbation on the inner product that the problem \eqref{Pblm Direct Formulation Perturbed} defines on $V$. The latter is reflected in the matrices \eqref{Eq matrix inner product} of the problem \eqref{Eq geometry effect on inner product problem} above.
\item
Although the operators $T(\zeta)$ converge in the strong operator topology as the perturbations $\zeta\in \pert(\Omega, \Gamma)$ converge, according to the hypothesis of theorem \ref{Th continuity of the conformal map operator}, it is an abstract statement. There are no explicit estimates depending on $\Vert \zeta\Vert_{\scriptscriptstyle W^{1, \infty}(\Gamma)}$, analogous to inequality \eqref{Ineq continuous dependence wrt left perturbation} presented in section \ref{Sec global estimate of perturbation}.
\item
The inner product that a geometric perturbation implicitly defines on the function space $V$ is the most important nonlinearity of the problem. It introduces a notion of orthogonality in the space which is hardly comparable with the standard one, beyond the topological equivalence of the induced norms. 
\item
The fractional mapping technique is a much simpler approach than the local charts strategy used in trace theorems. Its main contribution in this work, consists in exposing the challenges of the convergence rate question as well as the non-linearities involved, in a much neater way than the local charts approach.   
\item
The strong convergence statements attained in section \ref{Sec Strong Convergence}, realizing problem \eqref{Pblm Direct Formulation Original} as the strong limit of the family of problems \eqref{Pblm Direct Formulation Perturbed}, suggest numerical experimentation and a-posteriori estimates as the most feasible approach to gain insight in the convergence rate question as well as dependence with respect to the norms $W^{1, \infty}(\Gamma)$ and $C(\Gamma)$.
\item
The solution of the one dimensional case presented in section \ref{Sec One dimensional example}, using orthogonal decompositions on carefully chosen subspaces, illustrates the complexity of the convergence rate problem. This fact becomes even more dramatic due to the strong dependence on the one dimensional setting.
\item
In both settings, multiple and one dimensional, the necessity of testing the variational statements with functions of the structure \eqref{Def piecewise linear combination} (i.e. discontinuous across the interfaces), to obtain explicit a-priori estimates of the difference $p - q$ is self-evident. Such functions break the linear structure of the domain $V$ and obey to the disagreement of scaling coefficients in problems \eqref{Pblm Direct Formulation Original}, \eqref{Pblm Direct Formulation Perturbed}.
\item
Both of the mixed variational formulations: $\mathbf{L}^{\!2}$-$H^{1}$ and $\Hdiv$-$L^{\!2}$ demand coupling conditions on the function spaces for the trace on the interfaces $\Gamma, \Gammaz$. These conditions do not allow testing with discontinuous functions such as \ref{Def piecewise linear combination}. However, a mixed formulation setting $\mathbf{L}^{\!2}$-$H^{1}$ in one region, namely $\Omega_{1}$ and $\Hdiv$-$L^{\!2}$ in the other region, as the one introduced in \cite{MoralesShow2}, does not require continuity of the test functions across the interfaces. Hence, this is the formulation where convergence rate estimates (implicit or explicit, given the nonlinearity of the problem) can most likely be attained.
\end{asparaenum}
%
%
\section{Acknowledgments}
The author thanks to Universidad Nacional de Colombia, Sede Medell\'in for partially supporting this work under project HERMES 17194 and the Department of Energy USA for partially supporting this work by grant 98089.

\end{document}